\documentclass[a4paper,12pt]{amsart}
\usepackage[a4paper,hscale=0.7,vscale=0.75,centering]{geometry}
\usepackage{fullpage}
\usepackage{url}

\newtheorem{theorem}{Theorem}[section]
\newtheorem{lemma}[theorem]{Lemma}
\newtheorem{corollary}[theorem]{Corollary}

\theoremstyle{definition}

\newtheorem{example}[theorem]{Example}

\newtheorem{assumption}{Assumption}
\newtheorem*{proof0}{Proof of Theorem \ref{couplingTheorem}}
\newtheorem*{proof1}{Proof of Corollary \ref{corollary1}}
\newtheorem*{proof2}{Proof of Corollary \ref{corollary2}}
\newtheorem*{proof3}{Proof of Corollary \ref{corollary3}}

\theoremstyle{remark}
\newtheorem{remark}[theorem]{Remark}

\numberwithin{equation}{section}

\begin{document}

\title{Coupling and exponential ergodicity for stochastic differential equations driven by L\'{e}vy processes}

\author{Mateusz B. Majka}
\address{Institute for Applied Mathematics, University of Bonn, Endenicher Allee 60, 53115 Bonn, Germany}

\email{majka@uni-bonn.de}



\date{}

\dedicatory{}

\begin{abstract}
We present a novel idea for a coupling of solutions of stochastic differential equations driven by L\'{e}vy noise, inspired by some results from the optimal transportation theory. Then we use this coupling to obtain exponential contractivity of the semigroups associated with these solutions with respect to an appropriately chosen Kantorovich distance. As a corollary, we obtain exponential convergence rates in the total variation and standard $L^1$-Wasserstein distances.
\end{abstract}

\maketitle

\section{Introduction}

We consider stochastic differential equations of the form
\begin{equation}\label{SDE1}
 dX_t = b(X_t)dt + dL_t \,,
\end{equation}
where $(L_t)_{t \geq 0}$ is an $\mathbb{R}^d$-valued L\'{e}vy process and $b: \mathbb{R}^d \to \mathbb{R}^d$ is a continuous vector field satisfying a one-sided Lipschitz condition, i.e., there exists a constant $C_L > 0$ such that for all $x$, $y \in \mathbb{R}^d$ we have
\begin{equation}\label{onesidedLipschitz}
 \langle b(x) - b(y) , x - y \rangle \leq C_L |x-y|^2 \,.
\end{equation}
These assumptions are sufficient in order for (\ref{SDE1}) to have a unique strong solution (see Theorem 2 in \cite{gyongykrylov}). For any $t \geq 0$, denote the distribution of the random variable $L_t$ by $\mu_t$. Its Fourier transform $\widehat{\mu}_t$ is of the form
\begin{equation*}
\widehat{\mu}_t(z) = e^{t\psi(z)} \,, \; z \in \mathbb{R}^d \,,
\end{equation*}
where the \emph{L\'{e}vy symbol} (or \emph{L\'{e}vy exponent}) $\psi: \mathbb{R}^d \to \mathbb{C}$ is given by the L\'{e}vy - Khintchine formula (see e.g. \cite{apple} or \cite{sato}),
\begin{equation*}
\psi(z) = i\langle l, z \rangle - \frac{1}{2}\langle z, Az\rangle + \int_{\mathbb{R}^d} (e^{i\langle z, x \rangle} - 1 - i\langle z, x\rangle \mathbf{1}_{\{|x|\leq 1 \}})\nu(dx) \,,
\end{equation*}
for $z \in \mathbb{R}^d$. Here $l$ is a vector in $\mathbb{R}^d$, $A$ is a symmetric nonnegative-definite $d \times d$ matrix and $\nu$ is a measure on $\mathbb{R}^d$ satisfying
\begin{equation*}
\nu(\{0\}) = 0 \text{  and  } \int_{\mathbb{R}^d}(|x|^2 \wedge 1)\nu(dx) < \infty \,.
\end{equation*}
We call $(l, A, \nu)$ the \emph{generating triplet} of the L\'{e}vy process $(L_t)_{t\geq 0}$, whereas $A$ and $\nu$~are called, respectively, the \emph{Gaussian covariance matrix} and the \emph{L\'{e}vy measure} (or \emph{jump measure}) of $(L_t)_{t\geq 0}$.

In this paper we will be working with pure jump L\'{e}vy processes. We assume that in the generating triplet of $(L_t)_{t \geq 0}$ we have $l = 0$ and $A = 0$. By the L\'{e}vy - It\^{o} decomposition we know that there exists a Poisson random measure $N$ associated with $(L_t)_{t \geq 0}$ in such a way that
\begin{equation}\label{LevyViaPoisson}
  L_t = \int_0^t \int_{\{|v|>1\}} v N(ds,dv) + \int_0^t \int_{\{|v| \leq 1\}} v \widetilde{N}(ds,dv)\,,
\end{equation}
where
\begin{equation*}
\widetilde{N}(ds,dv) = N(ds,dv) - ds \, \nu(dv)
\end{equation*}
is the compensated Poisson random measure.

We will be considering the class of Kantorovich ($L^1$-Wasserstein) distances. For $p \geq 1$, we can define the $L^p$-\emph{Wasserstein distance} between two probability measures $\mu_1$ and $\mu_2$ on $\mathbb{R}^d$ by the formula
\begin{equation*}
 W_p(\mu_1, \mu_2) := \left( \inf_{\pi \in \Pi(\mu_1, \mu_2)} \int_{\mathbb{R}^d \times \mathbb{R}^d} \rho(x,y)^p \pi(dx \, dy) \right)^{\frac{1}{p}} \,,
\end{equation*}
where $\rho$ is a metric on $\mathbb{R}^d$ and $\Pi(\mu_1, \mu_2)$ is the family of all couplings of $\mu_1$ and $\mu_2$, i.e., $\pi \in \Pi(\mu_1, \mu_2)$ if and only if $\pi$ is a measure on $\mathbb{R}^{2d}$ having $\mu_1$ and $\mu_2$ as its marginals. We will be interested in the particular case of $p=1$ and the distance $\rho$ being given by a concave function $f: [0,\infty) \to [0,\infty)$ with $f(0) = 0$ and $f(x) > 0$ for $x > 0$ as
\begin{equation*}
 \rho(x,y) := f(|x-y|) \text{  for all } x, y \in \mathbb{R}^d \,.
\end{equation*}
We will denote the $L^1$-Wasserstein distance associated with a function $f$ by $W_f$. The most well-known examples are given by $f(x) = \mathbf{1}_{(0,\infty)}(x)$, which leads to the total variation distance (with $W_f(\mu_1, \mu_2) = \frac{1}{2}\| \mu_1 - \mu_2 \|_{TV}$) and by $f(x) = x$, which defines the standard $L^1$-Wasserstein distance (denoted later by $W_1$). For a detailed exposition of Wasserstein distances, see e.g. Chapter 6 in \cite{villani}.

For an $\mathbb{R}^d$-valued Markov process $(X_t)_{t \geq 0}$ with transition kernels $(p_t(x,\cdot))_{t \geq 0, x \in \mathbb{R}^d}$ we say that an $\mathbb{R}^{2d}$-valued process $(X_t',X_t'')_{t \geq 0}$ is a \emph{coupling} of two copies of the Markov process $(X_t)_{t \geq 0}$ if both $(X_t')_{t \geq 0}$ and $(X_t'')_{t \geq 0}$ are Markov processes with transition kernels $p_t$ but possibly with different initial distributions. We define the \emph{coupling time} $T$ for the marginal processes $(X_t')_{t \geq 0}$ and $(X_t'')_{t \geq 0}$ by $T := \inf \{ t \geq 0 : X_t' = X_t'' \}$. The coupling is called \emph{successful} if $T$~is almost surely finite. It is known (see e.g. \cite{lindvallbook} or \cite{thorisson}) that the condition 
\begin{equation*}
 \| \mu_1 p_t - \mu_2 p_t\|_{TV} \to 0 \text{  as } t \to \infty \text{  for any probability measures } \mu_1 \text{ and } \mu_2 \text{ on } \mathbb{R}^d
\end{equation*}
is equivalent to the property that for any two probability measures $\mu_1$ and $\mu_2$ on $\mathbb{R}^d$ there exist marginal processes $(X_t')_{t \geq 0}$ and $(X_t'')_{t \geq 0}$ with $\mu_1$ and $\mu_2$ as their initial distributions such that the coupling $(X_t',X_t'')_{t \geq 0}$ is successful. Here $\mu p_t(dy) = \int \mu(dx) p_t(x,dy)$.

Couplings of L\'{e}vy processes and related bounds in the total variation distance have recently attracted considerable attention. See e.g. \cite{bottcher}, \cite{schillingsztonyk} and \cite{schilling1} for couplings of pure jump L\'{e}vy processes, \cite{schilling2}, \cite{fywang} and \cite{jwang1} for the case of L\'{e}vy-driven Ornstein-Uhlenbeck processes and \cite{linwang}, \cite{jwang2} and \cite{song} for more general L\'{e}vy-driven SDEs with non-linear drift. See also \cite{kulik} and \cite{priola} for general considerations concerning ergodicity of SDEs with jumps. Furthermore, in a recent paper \cite{jwang4}, J. Wang investigated the topic of using couplings for obtaining bounds in the $L^p$-Wasserstein distances. 

Previous attempts at constructing couplings of L\'{e}vy processes or couplings of solutions to L\'{e}vy-driven SDEs include e.g. a coupling of subordinate Brownian motions by making use of the coupling of Brownian motions by reflection (see \cite{bottcher}), a coupling of compound Poisson processes obtained from certain couplings of random walks (see \cite{schilling1} for the original construction and \cite{jwang2} for a related idea applied to L\'{e}vy-driven SDEs) and a combination of the coupling by reflection and the synchronous coupling defined via its generator for solutions to SDEs driven by L\'{e}vy processes with a symmetric $\alpha$-stable component (see \cite{jwang4}). In the present paper we use a different idea for a coupling, as well as a different method of construction. Namely, we define a coupling by reflection modified in such a way that it allows for a positive probability of bringing the marginal processes to the same point if the distance between them is small enough. Such a behaviour makes it possible to obtain better convergence rates than a regular coupling by reflection, since it significantly decreases the probability that the marginal processes suddenly jump far apart once they have already been close to each other. We construct our coupling as a~solution to an explicitly given SDE, much in the vein of the seminal paper \cite{lindvall} by Lindvall and Rogers, where they constructed a coupling by reflection for diffusions with a drift. The formulas for the SDEs defining the marginal processes in our coupling are given by (\ref{eqX}) and (\ref{optimalY}) and the way we obtain them is explained in detail in Subsection \ref{subsectionConstruction}. Then, using this coupling, we construct a carefully chosen Kantorovich distance $W_f$~for an appropriate concave function $f$ such that
\begin{equation*}
 W_f(\mu_1 p_t , \mu_2 p_t) \leq e^{-ct} W_f(\mu_1, \mu_2)
\end{equation*}
holds for some constant $c > 0$ and all $t \geq 0$, where $\mu_1$ and $\mu_2$ are arbitrary probability measures on $\mathbb{R}^d$ and $(p_t)_{t \geq 0}$ is the transition semigroup associated with $(X_t)_{t \geq 0}$. Here $f$~and $c$ are mutually dependent and are chosen with the aim to make $c$ as large as possible, which leads to bounds that are in some cases close to optimal. A similar approach has been recently taken by Eberle in \cite{eberle}, where he used a specially constructed distance in order to investigate exponential ergodicity of diffusions with a drift. Historically, related ideas have been used e.g. by Chen and Wang in \cite{chen} and by Hairer and Mattingly in \cite{hairer}, to investigate spectral gaps for diffusion operators on $\mathbb{R}^d$ and to investigate ergodicity in infinite dimensions, respectively. It is important to point out that the distance function we choose is discontinuous. It is in fact of the form
\begin{equation*}
 f = f_1 + a\mathbf{1}_{(0,\infty)} \,,
\end{equation*}
where $f_1$ is a concave, strictly increasing $\mathcal{C}^2$ function with $f_1(0) = 0$, which from some point $R_1 > 0$ is extended in an affine way and $a$ is a positive constant. This choice of the distance (which is directly tied to our choice of the coupling) has an advantage in that it gives us upper bounds in both the total variation and standard $L^1$-Wasserstein distances (see Corollaries \ref{corollary1} and \ref{corollary2} and the discussion in Remark \ref{remarkDiscontinuityDiscussion}). 

Let us now state the assumptions that we will impose on the L\'{e}vy measure $\nu$ of the process $(L_t)_{t \geq 0}$.

\begin{assumption}\label{Assumption1}
 $\nu$ is rotationally invariant, i.e.,
 \begin{equation*}
 \nu(AB) = \nu(B)
\end{equation*}
for every Borel set $B \in \mathcal{B}(\mathbb{R}^d)$ and every $d \times d$ orthogonal matrix $A$.
\end{assumption}

\begin{assumption}\label{Assumption2}
 $\nu$ is absolutely continuous with respect to the Lebesgue measure on $\mathbb{R}^d$, with a density $q$ that is almost everywhere continuous on $\mathbb{R}^d$.
\end{assumption}

\begin{assumption}\label{Assumption3}
 There exist constants $m$, $\delta > 0$ such that $\delta < 2m$ and
 \begin{equation}\label{overlapCondition}
  \inf_{x \in \mathbb{R}^d : 0 < |x| \leq \delta} \int_{\{ |v| \leq m \} \cap \{ |v + x| \leq m \}} q(v) \wedge q(v + x) dv > 0 \,.
 \end{equation}
\end{assumption}

\begin{assumption}\label{Assumption4}
 There exists a constant $\varepsilon > 0$ such that $\varepsilon \leq \delta$ (with $\delta$ defined via (\ref{overlapCondition}) above) and
 \begin{equation*}
  \int_{\{ |v| \leq \varepsilon / 2 \}} q(v) dv > 0 \,.
 \end{equation*}
\end{assumption}

Assumptions \ref{Assumption1} and \ref{Assumption2} are used in the proof of Theorem \ref{couplingTheorem} to show that the solution to the SDE that we construct there is actually a coupling. Assumption \ref{Assumption1} is quite natural since we want to use reflection of the jumps. It is possible to extend our results to the case where the L\'{e}vy measure is only required to have a rotationally invariant component, but we do not do this in the present paper. Assumption \ref{Assumption3} is used in our calculations regarding the Wasserstein distances and is basically an assumption about sufficient overlap of the L\'{e}vy density $q$ and its translation. A related condition is used e.g. in \cite{schilling1} (see (1.3) in Theorem 1.1 therein) and in \cite{jwang1} to ensure that there is enough jump activity to provide a successful coupling. The restriction in (\ref{overlapCondition}) to the jumps bounded by $m$ is related to our coupling construction, see the discussion in Section \ref{subsectionConstruction}. Assumption \ref{Assumption4} ensures that we have enough small jumps to make use of the reflected jumps in our coupling (cf. the proof of Lemma \ref{lemma33}). All the assumptions together are satisfied by a large class of rotationally invariant L\'{e}vy processes, with symmetric $\alpha$-stable processes for $\alpha \in (0,2)$ being one of the most important examples. Note however, that our framework covers also the case of finite L\'{e}vy measures and even some cases of L\'{e}vy measures with supports separated from zero (see Example \ref{exampleAssumptions} for further discussion).

We must also impose some conditions on the drift function $b$. We have already assumed that it satisfies a one-sided Lipschitz condition, which guarantees the existence and uniqueness of a strong solution to (\ref{SDE1}). Now we define the function $\kappa: \mathbb{R}_{+} \to \mathbb{R}$ by setting $\kappa (|x-y|)$ to be the largest quantity such that
\begin{equation*}
 \langle b(x) - b(y) , x - y \rangle \leq -\kappa(|x-y|)|x-y|^2 \text{ for any } x, y \in \mathbb{R}^d \,,
\end{equation*}
and therefore it has to be defined as
\begin{equation}\label{defKappa}
 \kappa(r):=\inf \left\{ -\frac{\langle b(x) - b(y) , x - y \rangle}{|x-y|^2} : x,y \in \mathbb{R}^d \text{ such that } |x-y|=r \right\} \,.
\end{equation}
We have the following assumption.
\begin{assumption}\label{Assumption5}
 $\kappa$ is a continuous function satisfying
 \begin{equation*}
  \liminf_{r \to \infty} \kappa(r) > 0 \,.
 \end{equation*}
\end{assumption}
The above condition means that there exist constants $M > 0$ and $R > 0$ such that for all $x$, $y \in \mathbb{R}^d$ with $|x - y| \geq R$ we have
\begin{equation}\label{convexDrift}
 \langle b(x) - b(y) , x - y \rangle \leq -M|x-y|^2 \,.
\end{equation}
In other words, the drift $b$ is dissipative outside some ball of radius $R$. Note that if the drift is dissipative everywhere, i.e., (\ref{convexDrift}) holds for all $x$, $y \in \mathbb{R}^d$, then the proof of exponential convergence in the $L^1$-Wasserstein distance is quite straightforward, using just the synchronous coupling for $(L_t)_{t \geq 0}$ and the Gronwall inequality. Thus it is an interesting problem to try to obtain exponential convergence under some weaker assumptions on the drift.

We finally formulate our main results.

\begin{theorem}\label{couplingTheorem}
  Let us consider a stochastic differential equation
 \begin{equation}\label{SDEinTheorem}
  dX_t = b(X_t)dt + dL_t \,,
 \end{equation}
 where $(L_t)_{t \geq 0}$ is a pure jump L\'{e}vy process with the L\'{e}vy measure $\nu$ satisfying Assumptions \ref{Assumption1} and \ref{Assumption2}, whereas $b: \mathbb{R}^d \to \mathbb{R}^d$ is a continuous, one-sided Lipschitz vector field. Then a coupling $(X_t, Y_t)_{t \geq 0}$ of solutions to (\ref{SDEinTheorem}) can be constructed as a strong solution to the $2d$-dimensional SDE given by (\ref{eqX}) and (\ref{optimalY}), driven by a $d$-dimensional noise. If we additionally require Assumptions \ref{Assumption3}-\ref{Assumption5} to hold, then there exist a concave function $f$ and a constant $c > 0$ such that for any $t \geq 0$ we have
 \begin{equation}\label{mainInequality}
  \mathbb{E}f(|X_t - Y_t|) \leq e^{-ct} \mathbb{E}f(|X_0 - Y_0|)
 \end{equation}
and the coupling $(X_t, Y_t)_{t \geq 0}$ is successful.
\end{theorem}

Since the inequality (\ref{mainInequality}) holds for all couplings of the laws of $X_0$ and $Y_0$, directly from the definition of the Wasserstein distance $W_f$ we obtain the following result.

\begin{corollary}\label{mainTheorem}
 Let $(X_t)_{t \geq 0}$ be a solution to the SDE (\ref{SDEinTheorem}) with $(L_t)_{t \geq 0}$ and $b$ as in Theorem \ref{couplingTheorem}, satisfying Assumptions \ref{Assumption1}-\ref{Assumption5}. Then there exist a concave function $f$ and a constant $c > 0$ such that for any $t \geq 0$ and any probability measures $\mu_1$ and $\mu_2$ on $\mathbb{R}^d$ we have
\begin{equation}\label{Wfbound}
W_f(\mu_1 p_t, \mu_2 p_t) \leq e^{-ct}W_f(\mu_1, \mu_2) \,,
\end{equation}
where $(p_t)_{t \geq 0}$ is the semigroup associated with $(X_t)_{t \geq 0}$. 
\end{corollary}

The function $f$ in the theorem and the corollary above is given as $f = a \mathbf{1}_{(0,\infty)} + f_1$, where
\begin{equation}\label{fone}
 \begin{split}
  f_1(r) &= \int_0^r \phi(s)g(s) ds \\
  \phi(r) &= \exp{\left( -\int_0^r \frac{\bar{h}(t)}{C_{\varepsilon}} dt \right)} \,, \quad \bar{h}(r) = \sup_{t \in (r,r+\varepsilon)} t \kappa^{-}(t) \,,\\
  g(r) &= 1 - \frac{1}{2} \int_0^{r \wedge R_1} \frac{\Phi(t+\varepsilon)}{\phi(t)}dt \left( \int_0^{R_1} \frac{\Phi(t+\varepsilon)}{\phi(t)}dt \right)^{-1} \,, \quad \Phi(r) = \int_0^r \phi(s) ds \,,
 \end{split}
\end{equation}
while the contractivity constant $c$ is given by $c = \min \{ c_1/2K, \widetilde{C}_{\delta}/4 \}$ with
\begin{equation*}
c_1 = \frac{C_{\varepsilon}}{2}\left(\int_0^{R_1} \frac{\Phi(t+\varepsilon)}{\phi(t)} dt\right)^{-1} \quad \text{ and } \quad \widetilde{C}_{\delta} = \inf_{x \in \mathbb{R}^d : 0 < |x| \leq \delta} \int_{\mathbb{R}^d} q(v) \wedge q(v + x) dv \,.
\end{equation*}
Here $\kappa$ is the function defined by (\ref{defKappa}), the constants $R_0$ and $R_1$ are defined by
\begin{equation}\label{Rzero}
\begin{split}
 R_0 &= \inf \left\{R \geq 0 : \forall r \geq R : \kappa(r) \geq 0 \right\} \,, \\
 R_1 &= \inf \left\{ R \geq R_0 + \varepsilon : \forall r \geq R : \kappa(r) \geq \frac{2C_{\varepsilon}}{(R-R_0)R} \right\} \,,
 \end{split}
 \end{equation}
the constant $\delta$ comes from Assumption \ref{Assumption3}, the constant $\varepsilon \leq \delta$ comes from Assumption \ref{Assumption4} (see also Remark \ref{remarkChoiceEpsilonDelta}) and we have
\begin{equation}\label{Cepsilon}
 C_{\varepsilon} = 2 \int_{ -\varepsilon/4  }^0 |y|^2 \nu_1(dy) \,, \quad K =  \frac{C_L\delta + \widetilde{C}_{\delta}f_1(\delta)/2}{\widetilde{C}_{\delta}f_1(\delta)/2} \quad \text{ and } \quad a = Kf_1(\delta) \,,
\end{equation}
where $\nu_1$ is the first marginal of $\nu$ and the constant $C_L$ comes from (\ref{onesidedLipschitz}). Note that due to Assumptions \ref{Assumption3} and \ref{Assumption4} it is always possible to choose $\delta$ and $\varepsilon$ in such a way that $\widetilde{C}_{\delta} > 0$ and $C_{\varepsilon} > 0$ and due to Assumption \ref{Assumption5} the constants $R_0$ and $R_1$ are finite.

\begin{remark}
 The formulas for the function $f$ and the constant $c$ for which (\ref{Wfbound}) holds are quite sophisticated, but they are chosen in such a way as to try to make $c$ as large as possible and their choice is clearly motivated by the calculations in the proof, see Section \ref{Section3} for details. The contractivity constant $c$ can be seen to be in some sense close to optimal (at least in certain cases). See the discussion in Section \ref{Section4} for comparison of convergence rates in the $L^1$-Wasserstein distance in the case where the drift is assumed to be the gradient of a strongly convex potential and the case where convexity is only required to hold outside some ball.
\end{remark}

With the above notation and assumptions, we immediately get some important corollaries.

\begin{corollary}\label{corollary1}
For any $t \geq 0$ and any probability measures $\mu_1$ and $\mu_2$ on $\mathbb{R}^d$ we have
\begin{equation}\label{TVupperbound}
  \| \mu_1 p_t - \mu_2 p_t \|_{TV} \leq 2 a^{-1} e^{-ct} W_f(\mu_1, \mu_2) \,,
\end{equation}
where $a > 0$ is the constant defined by (\ref{Cepsilon}).
\end{corollary}

\begin{corollary}\label{corollary2}
For any $t \geq 0$ and any probability measures $\mu_1$ and $\mu_2$ on $\mathbb{R}^d$ we have
\begin{equation}\label{W1upperbound}
W_1(\mu_1 p_t , \mu_2 p_t) \leq 2\phi(R_0)^{-1} e^{-ct} W_f(\mu_1, \mu_2) \,,
\end{equation}
where the function $\phi$ and the constant $R_0 > 0$ are defined by (\ref{fone}) and (\ref{Rzero}), respectively.
\end{corollary}

\begin{remark}\label{remarkDiscontinuityDiscussion}
 The corollaries above follow in a straightforward way from (\ref{Wfbound}) by comparing the underlying distance function $f$ from below with the $\mathbf{1}_{(0,\infty)}$ function (corresponding to the total variation distance) and the identity function (corresponding to the standard $L^1$-Wasserstein distance), see Section \ref{Section4} for explicit proofs. In the paper \cite{eberle} by Eberle, which treated the diffusion case, a related concave function was constructed, although without a discontinuity at zero (and also extended in an affine way from some point). This leads to bounds of the form
 \begin{equation}\label{W1contraction}
  W_1(\mu_1 p_t, \mu_2 p_t) \leq L e^{-ct} W_1(\mu_1, \mu_2)
 \end{equation}
with some constants $L \geq 1$ and  $c > 0$, since such a continuous function $f$ can be compared with the identity function both from above and below. In our case we are not able to produce an inequality like (\ref{W1contraction}) due to the discontinuity at zero, but on the other hand we can obtain upper bounds (\ref{TVupperbound}) in the total variation distance, which is impossible in the framework of \cite{eberle}. Several months after the submission of the first version of the present manuscript, its author managed to modify the method presented here in order to obtain (\ref{Wfbound}) for L\'{e}vy-driven SDEs with a continuous function $f$ (which leads to (\ref{W1contraction})) by replacing Assumptions \ref{Assumption3} and \ref{Assumption4} with an assumption stating that the function $\varepsilon \mapsto \varepsilon/C_{\varepsilon}$ is bounded in a neighbourhood of zero (with $C_{\varepsilon}$ defined by (\ref{Cepsilon})), which is an assumption about sufficient concentration of the L\'{e}vy measure $\nu$ around zero (sufficient small jump activity, much higher than in the case of Assumptions \ref{Assumption3} and \ref{Assumption4}). This result was presented in \cite{majka16}, where trying to obtain the inequality (\ref{W1contraction}) was motivated by showing how it can lead to so-called $\alpha$-$W_1H$ transportation inequalities that characterize the concentration of measure phenomenon for solutions of SDEs of the form (\ref{SDE1}). The difference between the approach presented here and the approach in \cite{majka16} is in the method chosen to deal with the case in which the marginal processes in the coupling are already close to each other and contractivity can be spoilt by having undesirable large jumps. This can be dealt with either by introducing a discontinuity in the distance function and proceeding like in the proof of Lemma \ref{lemmaDiscontinuous} below or by making sure that we have enough small jumps. It is worth mentioning that in the meantime the inequality (\ref{W1contraction}) in the L\'{e}vy jump case was independently obtained by D. Luo and J. Wang in \cite{luowang}, by using a different coupling and under different assumptions (which are also, however, assumptions about sufficiently high small jump activity). In conclusion, it seems that in order to obtain (\ref{W1contraction}) one needs the noise to exhibit a diffusion-like type of behaviour (a lot of small jumps), while estimates of the type (\ref{TVupperbound}) and (\ref{W1upperbound}) can be obtained under much milder conditions.
\end{remark}

\begin{example}\label{exampleAssumptions}
 In order to better understand when Assumptions \ref{Assumption3} and \ref{Assumption4} are satisfied, let us examine a class of simple examples. We already mentioned that our assumptions hold for symmetric $\alpha$-stable processes with $\alpha \in (0,2)$, for which it is sufficient to take arbitrary $m > 0$ and arbitrary $\varepsilon = \delta < 2m$. Now let us consider one-dimensional L\'{e}vy measures of the form $\nu(dx) = \left( \mathbf{1}_{[-\theta,-\theta/\beta]}(x) + \mathbf{1}_{[\theta/\beta,\theta]}(x) \right) dx$ for arbitrary $\theta > 0$ and $\beta > 1$. If we would like the quantity appearing in Assumption \ref{Assumption3} to be positive, it is then best to take $m = \theta$. Note that if $\beta \leq 3$, then $2\theta/\beta \geq \theta - \theta/\beta$ (the gap in the support of $\nu$ is larger than the size of the part of the support contained in $\mathbb{R}_{+}$) and thus we need to have $\delta < \theta - \theta/\beta$ (taking $\delta = \theta - \theta/\beta$ or larger would result in an overlap of zero mass). This means that $\varepsilon/2 \leq \theta/2 - \theta/2\beta \leq \theta/\beta$ and thus the quantity in Assumption \ref{Assumption4} cannot be positive. On the other hand for $\beta > 3$ we can take any $\delta < 2\theta$ in Assumption \ref{Assumption3} and thus Assumption \ref{Assumption4} can also be satisfied.
 \end{example}

\begin{corollary}\label{corollary3}
 In addition to Assumptions \ref{Assumption1}-\ref{Assumption5}, suppose that the semigroup $(p_t)_{t \geq 0}$ preserves finite first moments, i.e., if a measure $\mu$ has a finite first moment, then for all $t > 0$ the measure $\mu p_t$ also has a finite first moment. Then there exists an invariant measure $\mu_{*}$ for the semigroup $(p_t)_{t \geq 0}$. Moreover, for any $t \geq 0$ and any probability measure $\eta$ we have
 \begin{equation}\label{Wfinvariantbound}
  W_f(\mu_{*}, \eta p_t) \leq e^{-ct}W_f(\mu_{*}, \eta)
 \end{equation}
and therefore
\begin{equation}\label{TVinvariantbound}
  \| \mu_{*} - \eta p_t \|_{TV} \leq 2 a^{-1} e^{-ct} W_f(\mu_{*}, \eta)
\end{equation}
and
\begin{equation}\label{W1invariantbound}
W_1(\mu_{*} , \eta p_t) \leq 2\phi(R_0)^{-1} e^{-ct} W_f(\mu_{*}, \eta) \,.
\end{equation}
\end{corollary}

To illustrate the usefulness of our approach, we can briefly compare our estimates with the ones obtained by other authors, who also investigated exponential convergence rates for semigroups $(p_t)_{t \geq 0}$ associated with solutions of equations like (\ref{SDE1}). In his recent paper \cite{song}, Y. Song obtained exponential upper bounds for $\|\delta_x p_t - \delta_y p_t \|_{TV}$ for $x$, $y \in \mathbb{R}^d$ using Malliavin calculus for jump processes, under some technical assumptions on the L\'{e}vy measure (which, however, does not have to be rotationally invariant) and under a global dissipativity condition on the drift. By our Corollary \ref{corollary1}, we get such bounds under a much weaker assumption on the drift. In \cite{jwang3}, J. Wang proved exponential ergodicity in the total variation distance for equations of the form (\ref{SDE1}) driven by $\alpha$-stable processes, while requiring the drift $b$ to satisfy a condition of the type $\langle b(x) , x \rangle \leq -C|x|^2$ when $|x| \geq R$ for some $R > 0$ and $C > 0$. In the proof he used a method involving the notions of $T$-processes and petite sets. His assumption on the drift is weaker than ours, but our results work for a much larger class of noise. Furthermore, in \cite{priola} the authors showed exponential ergodicity, again only in the $\alpha$-stable case, under some H\"{o}lder continuity assumptions on the drift, using two different approaches: by applying the Harris theorem and by a coupling argument. Kulik in \cite{kulik} also used a coupling argument to give some general conditions for exponential ergodicity, but in practice they can be difficult to verify. However, he gave a simple one-dimensional example of an equation like (\ref{SDE1}), with the drift satisfying a condition similar to the one in \cite{jwang3}, whose solution is exponentially ergodic under some relatively mild assumptions on the L\'{e}vy measure (see Proposition 0.1 in \cite{kulik}). It is important to point out that his results, similarly to ours, apply to some cases when the L\'{e}vy measure is finite (i.e., the equation (\ref{SDE1}) is driven by a compound Poisson process). All the papers mentioned above were concerned with bounds only in the total variation distance. On the other hand, J. Wang in \cite{jwang4} has recently obtained exponential convergence rates in the $L^p$-Wasserstein distances for the case when the noise in (\ref{SDE1}) has an $\alpha$-stable component and the drift is dissipative outside some ball. By our Corollary \ref{corollary2}, we get similar results in the $L^1$-Wasserstein distance for $\alpha$-stable processes with $\alpha \in (1,2)$, but also for a much larger class of L\'{e}vy processes without $\alpha$-stable components.

Several months after the previous version of the present manuscript had been submitted, a new paper \cite{luowang} by D. Luo and J. Wang appeared on arXiv. There the authors introduced yet another idea for a coupling of solutions to equations of the form (\ref{SDE1}) and used it to obtain exponential convergence rates for associated semigroups in both the total variation and the $L^1$-Wasserstein distances, as well as contractivity in the latter (cf. Remark \ref{remarkDiscontinuityDiscussion}). Their construction works under a technical assumption on the L\'{e}vy measure, which is essentially an assumption about its sufficient concentration around zero and it does not require the L\'{e}vy measure to be symmetric. However, the assumption in \cite{luowang} is significantly more restrictive than our Assumptions \ref{Assumption3} and \ref{Assumption4}. For example, it does not hold for finite L\'{e}vy measures as they do not have enough small jump activity, while our method works even in some cases where the support of the L\'{e}vy measure $\nu$ is separated from zero (cf. Example \ref{exampleAssumptions}).

The remaining part of this paper is organized as follows: In Section \ref{Section2} we explain the construction of our coupling and we formally prove that it is actually well defined. In Section \ref{Section3} we use it to prove the inequality (\ref{mainInequality}). In Section \ref{Section4} we prove Corollaries \ref{corollary1}, \ref{corollary2} and \ref{corollary3} and present some further calculations that provide additional insight into optimality of our choice of the contractivity constant $c$.

\section{Construction of the coupling}\label{Section2}

\subsection{Related ideas}

The idea for the coupling that we construct in this section comes from the paper \cite{mccann} by McCann, where he considered the optimal transport problem for concave costs on $\mathbb{R}$. Namely, given two probability measures $\mu_1$ and $\mu_2$ on $\mathbb{R}$, the problem is to find a measure $\gamma$ on $\mathbb{R}^2$ with marginals $\mu_1$ and $\mu_2$, such that the quantity
\begin{equation*}
 C(\gamma) := \int_{\mathbb{R}^2} c(x,y) d\gamma(x,y) \,,
\end{equation*}
called the transport cost, is minimized for a given concave function $c: \mathbb{R}^2 \to [0,\infty]$. McCann proved (see the remarks after the proof of Theorem 2.5 in \cite{mccann} and Proposition 2.12 therein) that the minimizing measure $\gamma$ (i.e., the optimal coupling of $\mu_1$ and $\mu_2$) is unique and independent of the choice of $c$, and gave an explicit expression for $\gamma$. Intuitively speaking, in the simplest case the idea behind the construction of $\gamma$ (i.e., of transporting the mass from $\mu_1$ to $\mu_2$) is to keep in place the common mass of $\mu_1$ and $\mu_2$ and to apply reflection to the remaining mass. McCann's paper only treats the one-dimensional case, but since in our setting the jump measure is rotationally invariant, it seems reasonable to try to use a similar idea for a coupling also in the multidimensional case. Note that we do not formally prove in this paper that the constructed coupling is in fact the optimal one. Statements like this are usually difficult to prove, but what we really need is just a good guess of how a coupling close to the optimal one should look. Then usefulness of the constructed coupling is verified by the good convergence rates that we obtain by its application.

A related idea appeared in the paper \cite{hsu} by Hsu and Sturm, where they dealt with couplings of Brownian motions, but the construction of what they call the mirror coupling can be also applied to other Markov processes. Assume we are given a symmetric transition density $p_t(x,z)$ on $\mathbb{R}$ and that we want to construct a coupling starting from $(x_1, x_2)$ as a joint distribution of an $\mathbb{R}^2$-valued random variable $\zeta = (\zeta_1, \zeta_2)$. We put
\begin{equation}\label{MirrorCoupling}
 \mathbb{P}(\zeta_2 = \zeta_1 | \zeta_1 = z_1) = \frac{p_t(x_1,z_1) \wedge p_t(x_2,z_1)}{p_t(x_1,z_1)}
\end{equation}
and
\begin{equation*}
 \mathbb{P}(\zeta_2 = x_1 + x_2 - \zeta_1 | \zeta_1 = z_1) = 1 - \frac{p_t(x_1,z_1) \wedge p_t(x_2,z_1)}{p_t(x_1,z_1)}
\end{equation*}
so the idea is that if the first marginal process moves from $x_1$ to $z_1$, then the second marginal can move either to the same point or to the point reflected with respect to $x_0 = \frac{x_1 + x_2}{2}$, with appropriately defined probabilities, taking into account the overlap of transition densities fixed at points $x_1$ and $x_2$. Alternatively, we can define this coupling by the joint transition kernel as
\begin{equation*}
 m_t(x_1,x_2,dy_1,dy_2) := \delta_{y_1}(dy_2)h_0(y_1)dy_1 + \delta_{Ry_1}(dy_2)h_1(y_1)dy_1 \,,
\end{equation*}
where $h_0(z) = p_t(x_1,z) \wedge p_t(x_2,z)$, $h_1(z) = p_t(x_1,z) - h_0(z)$ and $Ry_1 = x_1 + x_2 - y_1$. Hsu and Sturm prove that such a coupling is in fact optimal for concave, strictly increasing cost functions.

Now let us also recall the ideas from \cite{lindvall} by Lindvall and Rogers, where they constructed a coupling $(X_t,Y_t)_{t \geq 0}$ by reflection for diffusions by defining the second marginal process $(Y_t)_{t \geq 0}$ as a solution to an appropriate SDE. If we have a stochastic differential equation 
\begin{equation}\label{BrownianX}
 dX_t = b(X_t)dt + dB_t
\end{equation}
driven by a $d$-dimensional Brownian motion $(B_t)_{t \geq 0}$, we can define $(Y_t)_{t \geq 0}$ by setting
\begin{equation}\label{couplingSDE}
 dY_t = b(Y_t)dt + (I - 2e_{t}e_{t}^T)dB_t \,,
\end{equation}
where
\begin{equation}\label{defe}
 e_t := \frac{X_t - Y_t}{|X_t- Y_t|} \,.
\end{equation}
Of course, the equation (\ref{couplingSDE}) only makes sense for $t < T$, where $T:= \inf \{t \geq 0 : X_t = Y_t \}$, but we can set $Y_t := X_t$ for $t \geq T$. The proof that the equations (\ref{BrownianX}) and (\ref{couplingSDE}) together define a coupling, i.e., the solution $(Y_t)_{t \geq 0}$ to the equation (\ref{couplingSDE}) has the same finite dimensional distributions as the solution $(X_t)_{t \geq 0}$ to the equation (\ref{BrownianX}), is quite simple in the Brownian setting. It is sufficient to use the L\'{e}vy characterization theorem for Brownian motion, since the process $A_t := I - 2e_te_t^T$ takes values in orthogonal matrices (and thus the process $(\tilde{B}_t)_{t \geq 0}$ defined by $d\tilde{B}_t := A_tdB_t$ is also a Brownian motion).

Similarly, if we consider an equation like (\ref{couplingSDE}) but driven by a rotationally invariant L\'{e}vy process $(L_t)_{t \geq 0}$ instead of the Brownian motion, it is possible to show that the process $(\tilde{L}_t)_{t \geq 0}$ defined by $d\tilde{L}_t := A_{t-}dL_t$ with $A_{t-} := I - 2e_{t-}e_{t-}^T$ is a L\'{e}vy process with the same finite dimensional distributions as $(L_t)_{t \geq 0}$. However, a corresponding coupling by reflection for L\'{e}vy processes would not be optimal and we were not able to obtain contractivity in any distance $W_f$ using this coupling. Intuitively, this follows from the fact that such a construction allows for a situation in which two jumping processes, after they have already been close to each other, suddenly jump far apart. We need to somehow restrict such behaviour and therefore we use a more sophisticated construction.

\subsection{Construction of the SDE}\label{subsectionConstruction}

We apply the ideas from \cite{mccann} and \cite{hsu} by coupling the jumps of $(X_t)_{t \geq 0}$ and $(Y_t)_{t \geq 0}$ in an appropriate way. Namely, we would like to use the coupling by reflection modified in such a way that it allows for a positive probability of $(Y_t)_{t \geq 0}$ jumping to the same point as $(X_t)_{t \geq 0}$. In order to employ this additional feature, we need to modify the Poisson random measure $N$ associated with $(L_t)_{t \geq 0}$ via (\ref{LevyViaPoisson}). Recall that there exists a sequence $(\tau_j)_{j=1}^{\infty}$ of random variables in $\mathbb{R}_{+}$ encoding the jump times and a~sequence $(\xi_j)_{j=1}^{\infty}$ of random variables in $\mathbb{R}^d$ encoding the jump sizes such that
\begin{equation*}
 N((0,t],A)(\omega) = \sum_{j=1}^{\infty} \delta_{(\tau_j(\omega),\xi_j(\omega))}((0,t] \times A) \text{ for all } \omega \in \Omega \text{ and } A \in \mathcal{B}(\mathbb{R}^d) 
\end{equation*}
(see e.g. \cite{peszat}, Chapter 6). At the jump time $\tau_j$ the process $(X_t)_{t \geq 0}$ jumps from the point $X_{\tau_j-}$ to $X_{\tau_j}$ and our goal is to find a way to determine whether the jump of $(Y_t)_{t \geq 0}$ should be reflected or whether $(Y_t)_{t \geq 0}$ should be forced to jump to the same point that $(X_t)_{t \geq 0}$ jumped to. In order to achieve this, let us observe that instead of considering the Poisson random measure $N$ on $\mathbb{R}_{+} \times \mathbb{R}^d$, we can extend it to a Poisson random measure on $\mathbb{R}_{+} \times \mathbb{R}^d \times [0,1]$, replacing the $d$-dimensional random variables $\xi_j$ determining the jump sizes of $(L_t)_{t \geq 0}$, with the $(d+1)$-dimensional random variables $(\xi_j, \eta_j)$, where each $\eta_j$ is a uniformly distributed random variable on $[0,1]$. Thus we have
\begin{equation*}
  N((0,t],A)(\omega) = \sum_{j=1}^{\infty} \delta_{(\tau_j(\omega),\xi_j(\omega),\eta_j(\omega))}((0,t] \times A \times [0,1]) \text{ for all } \omega \in \Omega \text{ and } A \in \mathcal{B}(\mathbb{R}^d) 
\end{equation*}
and by a slight abuse of notation we can write
\begin{equation}\label{compensatedLevy}
 L_t = \int_0^t \int_{\{|v|>1\} \times [0,1]} v N(ds,dv,du) + \int_0^t \int_{\{|v| \leq 1\} \times [0,1]} v \widetilde{N}(ds,dv,du)\,,
\end{equation}
denoting our extended Poisson random measure also by $N$. With this notation, if there is a jump at time $t$, then the process $(X_t)_{t \geq 0}$ moves from the point $X_{t-}$ to $X_{t-} + v$ and we draw a random number $u \in [0,1]$ which is then used to determine whether the process $(Y_t)_{t \geq 0}$ should jump to the same point that $(X_t)_{t \geq 0}$ jumped to, or whether it should be reflected just like in the ``pure'' reflection coupling. In order to make this work, we introduce a control function $\rho$ with values in $[0,1]$ that will determine the probability of bringing the processes together. Our idea is based on the formula (\ref{MirrorCoupling}) and uses the minimum of the jump density $q$ and its translation by the difference of the positions of the two coupled processes before the jump time, that is, by the vector
\begin{equation*}
 Z_{t-} := X_{t-} - Y_{t-} \,.
\end{equation*}
Our first guess would be to define our control function by
\begin{equation}\label{defRho}
 \rho(v, Z_{t-}) := \min \left\{ \frac{q(v + Z_{t-})}{q(v)}, 1 \right\} = \frac{q(v + Z_{t-}) \wedge q(v)}{q(v)}
\end{equation}
when $q(v) > 0$. We set $\rho(v, Z_{t-}) := 1$ if $q(v) = 0$. Note that we have $q(v+Z_{t-})/q(v) = q(v+X_{t-}-Y_{t-})/q(v+X_{t-}-X_{t-})$, so we can look at this formula as comparing the translations of $q$ by the vectors $Y_{t-}$ and $X_{t-}$, respectively. The idea here is that ``on average'' the probability of bringing the processes together should be equal to the ratio of the overlapping mass of the jump density $q$ and its translation and the total mass of $q$. However, for technical reasons, we will slightly modify this definition. 

Namely, we will only apply our coupling construction presented above to the jumps of size bounded by a constant $m > 0$ satisfying Assumption \ref{Assumption3}. For the larger jumps we will apply the synchronous coupling, i.e., whenever $(X_t)_{t \geq 0}$ makes a jump of size greater than $m$, we will let $(Y_t)_{t \geq 0}$ make exactly the same jump. The rationale behind this is the following. First, this modification allows us to control the size of jumps of the difference process $Z_t := X_t - Y_t$. If $(X_t)_{t \geq 0}$ makes a large jump $v$, then instead of reflecting the jump for $(Y_t)_{t \geq 0}$ and having a large change in the value of $Z_t$, we make the same jump $v$ with $(Y_t)_{t \geq 0}$ and the value of $Z_t$ does not change at all. Secondly, by doing this we do not in any way spoil the contractivity in $W_f$ that we want to show. As will be evident in the proof, what is crucial for the contractivity is on one hand the reflection applied to small jumps only (see Lemma \ref{lemma33} and Lemma \ref{lemmaFunctionalInequality}) and on the other the quantity (\ref{overlapCondition}) from Assumption \ref{Assumption3} (see Lemma \ref{lemmaDiscontinuous}). If the latter, however, holds for some $m_0 > 0$ then it also holds for all $m \geq m_0$ and in our calculations we can always choose $m$ large enough if needed (see the inequality (\ref{choiceOfm1}) in the proof of Lemma \ref{lemma33} and (\ref{choiceOfm2}) after the proof of Lemma \ref{lemmaDiscontinuous}). Therefore choosing a large but finite $m$ is a better solution than constructing a coupling with $m = \infty$ (i.e., applying our ``mirror'' construction to jumps of all sizes), which would require us to impose an additional assumption on the size of jumps of the noise $(L_t)_{t \geq 0}$.

Now that we have justified making such an adjustment, note that for any fixed $m > 1$ we can always write (\ref{compensatedLevy}) as
\begin{equation*}
\begin{split}
 L_t &= \int_0^t \int_{\{ |v| > m \} \times [0,1]} v N(ds,dv,du) + \int_0^t \int_{\{ |v| \leq m \} \times [0,1]} v \widetilde{N}(ds,dv,du)\\
 &+ \int_0^t \int_{\{ m \geq |v| > 1 \} \times [0,1]} v \nu(dv) du ds \,.
 \end{split}
\end{equation*}
Then we can include the last term appearing above in the drift $b$ in the equation (\ref{SDE1}) describing $(X_t)_{t \geq 0}$. Obviously such a change of the drift does not influence its dissipativity properties. Thus, once we have fixed a large enough $m$ (see the discussion above), we can for notational convenience redefine $(L_t)_{t \geq 0}$ and $b$ by setting
\begin{equation}\label{newLt}
 L_t := \int_0^t \int_{\{ |v| > m \} \times [0,1]} v N(ds,dv,du) + \int_0^t \int_{\{ |v| \leq m \} \times [0,1]} v \widetilde{N}(ds,dv,du)
\end{equation}
and modifying $b$ accordingly.

Since we want to apply different couplings for the compensated and uncompensated parts of $(L_t)_{t \geq 0}$, we actually need to modify the definition (\ref{defRho}) of the control function $\rho$ by putting
\begin{equation*}
 \rho(v, Z_{t-}) := \frac{q(v) \wedge q(v + Z_{t-})\mathbf{1}_{\{ |v+Z_{t-}| \leq m \}}}{q(v)}
\end{equation*}
Observe that with our new definition for any integrable function $f$ and any $z \in \mathbb{R}^d$ we have
\begin{equation*}
 \begin{split}
  \int_{\{ |v| \leq m \}} f(v)\rho(v,z) \nu(dv) &= \int_{\{ |v| \leq m \}} f(v) \frac{q(v) \wedge q(v+z)\mathbf{1}_{\{ |v + z| \leq m \}}}{q(v)} q(v) dv \\
  &= \int_{\{ |v| \leq m \} \cap \{ |v + z| \leq m \}} f(v)\left(q(v) \wedge q(v+z)\right) dv \,,
 \end{split}
\end{equation*}
while with (\ref{defRho}) we would just have
\begin{equation*}
  \int_{\{ |v| \leq m \}} f(v)\rho(v,z) \nu(dv) = \int_{\{ |v| \leq m \}} f(v)\left(q(v) \wedge q(v+z)\right) dv \,.
\end{equation*}
We will use this fact later in the proof of Lemma \ref{charFunctionsTheorem}. On an intuitive level, if the distance $Z_{t-}$ between the processes before the jump is big (much larger than $m$), and we are only considering the jumps bounded by $m$ (and thus $|v+Z_{t-}|$ is still big), then the probability of bringing the processes together should be zero, while the quantity (\ref{defRho}) can still be positive in such a situation. The restriction we introduce in the definition of $\rho$ eliminates this problem.

To summarize, in our construction once we have the number $u \in [0,1]$, if the jump vector of $(X_t)_{t \geq 0}$ at time $t$ is $v$ and $|v| \leq m$, then the jump vector of $(Y_t)_{t \geq 0}$ should be $X_{t-} - Y_{t-} + v$ (so that $(Y_t)_{t \geq 0}$ jumps from $Y_{t-}$ to $X_{t-} + v$) when
\begin{equation}\label{uMniejsze}
 u < \rho(v, Z_{t-}) \,.
\end{equation}
Otherwise the jump of $(Y_t)_{t \geq 0}$ should be $v$ reflected with respect to the hyperplane spanned by the vector $e_{t-} = (X_{t-} - Y_{t-})/|X_{t-} - Y_{t-}|$. If $|v| > m$, then the jump of $(Y_t)_{t \geq 0}$ is the same as the one of $(X_t)_{t \geq 0}$, i.e., it is also given by the vector $v$.

We are now ready to define our coupling by choosing an appropriate SDE for the process $(Y_t)_{t \geq 0}$. 
Recall that $(X_t)_{t \geq 0}$ is given by (\ref{SDE1}) and thus
\begin{equation}\label{eqX}
 dX_t = b(X_t)dt + \int_{\{ |v| > m \} \times [0,1]} v N(dt,dv,du) + \int_{\{ |v| \leq m \} \times [0,1]} v \widetilde{N}(dt,dv,du) \,.
\end{equation}
Now, in view of the above discussion, we consider the SDE
\begin{equation}\label{optimalY}
 \begin{split}
  dY_t &= b(Y_t)dt + \int_{\{ |v| > m \} \times [0,1]} v N(dt,dv,du) \\
  &+ \int_{\{ |v| \leq m \} \times [0,1]} (X_{t-} - Y_{t-} + v)\mathbf{1}_{\{ u < \rho(v, Z_{t-}) \}} \widetilde{N}(dt,dv,du) \\
  &+ \int_{\{ |v| \leq m \} \times [0,1]} R(X_{t-},Y_{t-})v\mathbf{1}_{\{ u \geq \rho(v, Z_{t-}) \}} \widetilde{N}(dt,dv,du) \,,
 \end{split}
\end{equation}
where 
\begin{equation*}
 R(X_{t-},Y_{t-}) := I - 2\frac{(X_{t-}-Y_{t-})(X_{t-}-Y_{t-})^T}{|X_{t-} - Y_{t-}|^2} = I - 2e_{t-}e_{t-}^T
\end{equation*}
is the reflection operator like in (\ref{couplingSDE}) with $e_{t}$ defined by (\ref{defe}). Observe that if $Z_{t-} = 0$, then $\rho(v, Z_{t-}) = 1$ and the condition (\ref{uMniejsze}) is satisfied almost surely, so after $Z_t$ hits zero once, it stays there forever. Thus, if we denote
\begin{equation}\label{theCouplingTime}
 T := \inf \{ t \geq 0 : X_t = Y_t \} \,,
\end{equation}
then $X_t = Y_t$ for any $t \geq T$.

We can equivalently write (\ref{optimalY}) in a more convenient way as
\begin{equation}\label{optimalYequation}
  \begin{split}
  dY_t &= b(Y_t)dt + \int_{\{ |v| > m \} \times [0,1]} v N(dt,dv,du) \\
  &+ \int_{\{ |v| \leq m \} \times [0,1]} R(X_{t-},Y_{t-})v \widetilde{N}(dt,dv,du) \\
  &+ \int_{\{ |v| \leq m \} \times [0,1]} (X_{t-} - Y_{t-} + v - R(X_{t-},Y_{t-})v)\mathbf{1}_{\{ u < \rho(v, Z_{t-}) \}} \widetilde{N}(dt,dv,du) \,.
 \end{split}
\end{equation}

\subsection{Auxiliary estimates}

At first glance, it is not clear whether the above equation even has a solution or if $(X_t,Y_t)_{t \geq 0}$ indeed is a coupling. Before we answer these questions, we will first show some estimates of the coefficients of (\ref{optimalYequation}), which will be useful in the sequel (see Lemmas \ref{charFunctionsTheorem} and \ref{lemma32}).

\begin{lemma}\label{linearGrowthLemma}
(Linear growth) There exists a constant $C = C(m) > 0$ such that for any $x$, $y \in \mathbb{R}^d$ we have
 \begin{equation*}
  \int_{\{ |v| \leq m \} \times [0,1]} |x - y + v - R(x,y)v|^2 \mathbf{1}_{\{ u < \rho(v, x - y) \}} \nu(dv)du \leq C(1+|x - y|^2) \,.
 \end{equation*}
 \begin{proof}
  We will keep using the notation $z = x - y$. We have
  \begin{equation}\label{linGrowthEst}
  \begin{split}
   \int_{\{ |v| \leq m \}} |z + v - R(x,y)v|^2 \rho(v, z) \nu(dv) & \leq 2\int_{\{ |v| \leq m \}} |z + v|^2 \rho(v, z) \nu(dv)\\
   &+ 2\int_{\{ |v| \leq m \}} |R(x,y)v|^2 \rho(v, z) \nu(dv)
   \end{split}
  \end{equation}
and, since $R$ is an isometry, we can estimate
\begin{equation*}
  \begin{split}
   2\int_{\{ |v| \leq m \}} &|R(x,y)v|^2 \rho(v, z) \nu(dv) = 2\int_{\{ |v| \leq m \}} |v|^2 \rho(v, z) \nu(dv) \\
   &\leq 2\int_{\{ |v| \leq m \}} |v|^2 q(v + z) \wedge q(v) dv \leq 2\int_{\{ |v| \leq m \}} |v|^2 q(v) dv = 2\int_{\{ |v| \leq m \}} |v|^2 \nu(dv) \,.
   \end{split}
  \end{equation*}
  The last integral is of course finite, since $\nu$ is a L\'{e}vy measure. We still have to bound the first integral on the right hand side of (\ref{linGrowthEst}). We have
  \begin{equation*}
  \begin{split}
   2\int_{\{ |v| \leq m \}} |z + v|^2 \rho(v, z) \nu(dv) &\leq 2\int_{\{ |v| \leq m \}} |z + v|^2  q(v + z) \wedge q(v) dv \\
   &= 2\int_{\{ |v-z| \leq m \}} |v|^2  q(v) \wedge q(v - z) dv \,.
   \end{split}
  \end{equation*}
  Now let us consider two cases. First assume that $|z| \leq 2m$ (instead of $2$ we can also take any positive number strictly greater than $1$). Then
  \begin{equation*}
   2\int_{\{ |v-z| \leq m \}} |v|^2  q(v) \wedge q(v - z) dv \leq 2\int_{\{ |v-z| \leq m \}} |v|^2  \nu(dv) \leq 2\int_{\{ |v| \leq 3m \}} |v|^2  \nu(dv) < \infty \,.
  \end{equation*}
  On the other hand, when $|z| > 2m$, we have
  \begin{equation*}
   \{ v \in \mathbb{R}^d : |v-z| \leq m \} \subset \{ v \in \mathbb{R}^d : |v| \leq m \}^c =: B(m)^c \,,
  \end{equation*}
  and $\nu(B(m)^c) < \infty$, which allows us to estimate
  \begin{equation*}
   \begin{split}
    2&\int_{\{ |v-z| \leq m \}} |v|^2  q(v) \wedge q(v - z) dv \\
    &\leq 4\int_{\{ |v-z| \leq m \}} |v-z|^2  q(v) \wedge q(v - z) dv + 4\int_{\{ |v-z| \leq m \}} |z|^2  q(v) \wedge q(v - z) dv\\
    &\leq 4\int_{\{ |v-z| \leq m \}} |v-z|^2 q(v - z) dv + 4\int_{\{ |v-z| \leq m \}} |z|^2  q(v) dv\\
    &\leq 4\int_{\{ |v| \leq m \}} |v|^2 \nu(dv) + 4|z|^2 \nu(B(m)^c) \,.
   \end{split}
  \end{equation*}
  Hence, by choosing
  \begin{equation*}
   C:=\max \left\{ 2\int_{\{ |v| \leq 3m \}} |v|^2 \nu(dv) + 2\int_{\{ |v| \leq m \}} |v|^2 \nu(dv) , 6\int_{\{ |v| \leq m \}} |v|^2 \nu(dv), 4\nu(B(m)^c)\right\}
  \end{equation*}
  we get the desired result.
 \end{proof}
\end{lemma}
Here we should remark that by the above lemma we have
\begin{equation*}
  \mathbb{P}\left( \int_0^t \int_{\{ |v| \leq m \} \times [0,1]} |Z_{s-} + v - R(X_{s-},Y_{s-})v|^2 \mathbf{1}_{\{ u < \rho(v, Z_{s-}) \}} \nu(dv)du ds < \infty \right) = 1 \,.
\end{equation*}
We will use this fact later on.

The next thing we need to show is that the (integrated) coefficients are continuous in the solution variable. Note that obviously
\begin{equation*}
 \int_{\{ |v| \leq m \} \times [0,1]} |R(x+h,y)v - R(x,y)v|^2 \nu(dv)du \to 0 \,, \text{ as } h \to 0 \,,
\end{equation*}
so we just need to take care of the part involving $\rho(v,z)$. Before we proceed though, let us make note of the following fact.
\begin{remark}\label{RemarkAboutRho}
 For a fixed value of $z \neq 0$, the measure
 \begin{equation*}
 \rho(v,z)\nu(dv)
 \end{equation*}
 is a finite measure on $\mathbb{R}^d$. Indeed, if $z \neq 0$, we can choose a neighbourhood $U$ of $z$ such that $0 \notin \bar{U}$. Then $U - z$ is a neighbourhood of $0$ and we have
 \begin{equation*}
  \begin{split}
   \int_{\mathbb{R}^d} \rho(v,z)\nu(dv) &= \int_{U} \rho(v,z)\nu(dv) + \int_{U^c} \rho(v,z)\nu(dv) \\
   &\leq \int_{U} q(v)dv + \int_{U^c} q(v+z)dv \\
   &= \int_{U} q(v)dv + \int_{(U - z)^c} q(v)dv < \infty \,,
  \end{split}
 \end{equation*}
since $\nu$ is a L\'{e}vy measure.
\end{remark}

\begin{lemma}\label{continuityLemma}
(Continuity condition) For any $x$, $y \in \mathbb{R}^d$ and $z = x - y$ we have
 \begin{equation*}
 \begin{split}
 \int_{\{ |v| \leq m \} \times [0,1]} |(x &+ h - y + v - R(x+h,y)v)\mathbf{1}_{\{ u < \rho(v,z+h) \}} \\
 &- (x - y + v - R(x,y)v)\mathbf{1}_{\{ u < \rho(v,z) \}}|^2 \nu(dv)du \to 0 \,, \text{ as  } h \to 0 \,.
 \end{split}
 \end{equation*}
 \begin{proof}
  We have
 \begin{equation*}
 \begin{split}
 \int_{\{ |v| \leq m \} \times [0,1]} &|(x + h - y + v - R(x+h,y)v)\mathbf{1}_{\{ u < \rho(v,z+h) \}} \\
 &- (x - y + v - R(x,y)v)\mathbf{1}_{\{ u < \rho(v,z) \}}|^2 \nu(dv)du \\
 &= \int_{\{ |v| \leq m \} \times [0,1]} |(x + h - y + v - R(x+h,y)v)\mathbf{1}_{\{ u < \rho(v,z+h) \}} \\
 &- (x - y + v - R(x,y)v)\mathbf{1}_{\{ u < \rho(v,z+h) \}} \\
 &+ (x - y + v - R(x,y)v)\mathbf{1}_{\{ u < \rho(v,z+h) \}} \\
 &- (x - y + v - R(x,y)v)\mathbf{1}_{\{ u < \rho(v,z) \}}|^2 \nu(dv)du \\
 &\leq 2\int_{\{ |v| \leq m \}} |h - R(x+h,y)v + R(x,y)v|^2\rho(v,z+h) \nu(dv) \\
 &+ 2\int_{\{ |v| \leq m \}} |x - y + v - R(x,y)v|^2|\rho(v,z+h) - \rho(v,z)| \nu(dv)\\
 &=: I_1 + I_2 \,.
 \end{split}
 \end{equation*}
 Taking into account Remark \ref{RemarkAboutRho} and using the dominated convergence theorem, we can easily show that $I_1$ converges to zero when $h \to 0$. As for $I_2$, observe that
 \begin{equation*}
 \begin{split}
  |\rho(v,z+h)& - \rho(v,z)|\mathbf{1}_{\{ |v| \leq m \}} \\
  &= \frac{|q(v+z+h)\mathbf{1}_{\{ |v+z+h| \leq m \}} \wedge q(v) - q(v+z)\mathbf{1}_{\{ |v+z| \leq m \}} \wedge q(v)|}{|q(v)|} \mathbf{1}_{\{ |v| \leq m \}} \,.
  \end{split}
 \end{equation*}
Recall that by Assumption \ref{Assumption2}, the density $q$ is continuous almost everywhere on $\mathbb{R}^d$. Moreover, for a fixed $z \in \mathbb{R}^d$ the function $\mathbf{1}_{\{ |v+z| \leq m \}}$ is continuous outside of the set $\{v \in \mathbb{R}^d : |v+z| = m \}$, which is of measure zero. Therefore, using the dominated convergence theorem once again, we show that $I_2 \to 0$ when $h \to 0$.
 \end{proof}
\end{lemma}

\subsection{Existence of a solution}\label{sectionExistence}

Note that having the above estimates, it would be possible to prove existence of a weak solution to the $2d$-dimensional system given by (\ref{eqX}) and (\ref{optimalY}), using Theorem 175 in \cite{siturong}. However, there is a simpler method allowing to prove even more, namely, existence of a unique strong solution. To this end, we will use the so-called interlacing technique. This technique of modifying the paths of a process by adding jumps defined by a Poisson random measure of finite intensity is well known, cf. e.g. Theorem IV-9.1 in \cite{ikeda} or Theorem 6.2.9 in \cite{apple}. We first notice that without loss of generality it allows us to focus on the small jumps of size bounded by $m$, as we can always add the big jumps later, both to $(X_t)_{t \geq 0}$ and $(Y_t)_{t \geq 0}$. Hence we can consider the equation for $(Y_t)_{t \geq 0}$ written as
\begin{equation}\label{optimalYequation2}
   \begin{split}
 dY_t &= b(Y_t)dt + \int_{\{ |v| \leq m \} \times [0,1]} R(X_{t-},Y_{t-})v \widetilde{N}(dt,dv,du) \\
  &+ \int_{\{ |v| \leq m \} \times [0,1]} (X_{t-} - Y_{t-} + v - R(X_{t-},Y_{t-})v)\mathbf{1}_{\{ u < \rho(v, Z_{t-}) \}} \widetilde{N}(dt,dv,du) \,.
 \end{split}
\end{equation}
Now observe that if we only consider the equation
\begin{equation}\label{YwithoutJumpIntegral}
 dY_t^1 = b(Y_t^1)dt + \int_{\{ |v| \leq m \} \times [0,1]} R(X_{t-},Y_{t-}^1)v \widetilde{N}(dt,dv,du) \,,
\end{equation}
it is easy to see that it has a unique strong solution since the process $(X_t,Y_t^1)_{t \geq 0}$ up to its coupling time $T$ takes values in the region of $\mathbb{R}^{2d}$ in which the function $R$ is locally Lipschitz and has linear growth. Then note that the second integral appearing in (\ref{optimalYequation2}) represents a sum of jumps of which (almost surely) there is only a finite number on any finite time interval, since
\begin{equation*}
 \int_{\mathbb{R}^d \times [0,1]} \mathbf{1}_{\{ u < \rho(v, Z_{t-}) \}} \nu(dv) du = \int_{\mathbb{R}^d} \rho(v, Z_{t-}) \nu (dv) < \infty \,,
\end{equation*}
as long as $Z_{t-} \neq 0$ (see Remark \ref{RemarkAboutRho} above). Then in principle in such situations it is possible to use the interlacing technique to modify the paths of the process $(Y_t^1)_{t \geq 0}$ by adding the jumps defined by the second integral in (\ref{optimalYequation2}), see e.g. the proof of Proposition 2.2 in \cite{luowang} for a similar construction. Here, however, our particular case is even simpler. Namely, let us consider a uniformly distributed random variable $\xi \in [0,1]$ and define
\begin{equation*}
 \tau_1 := \inf \{ t > 0 : \xi < \rho(\Delta L_t , Z^1_{t-}) \} \,,
\end{equation*}
where $Z^1_{t} : = X_{t} - Y_{t}^1$ and $(L_t)_{t \geq 0}$ is the L\'{e}vy process associated with $N$. Then if we define a process $(Y^2_t)_{t \geq 0}$ by adding the jump of size 
$X_{{\tau_{1}}-} - Y_{{\tau_{1}}-}^1 + \Delta L_{\tau_{1}} - R(X_{{\tau_{1}}-},Y_{{\tau_{1}}-}^1)\Delta L_{\tau_{1}}$ to the path of $(Y_t^1)_{t \geq 0}$ at time $\tau_1$, we see that $Y^2_{\tau_1} = X_{\tau_1}$. Moreover, since $\rho(v,0) = 1$ for any $v \in \mathbb{R}^d$, we have $Y^2_t = X_t$ for all $t \geq \tau_1$. Thus we only need to add one jump to the solution of (\ref{YwithoutJumpIntegral}) in order to obtain a process which behaves like a solution to (\ref{optimalYequation2}) up to the coupling time, and like the process $(X_t)_{t \geq 0}$ later on. In consequence we obtain a solution $(X_t,Y_t)_{t \geq 0}$ to the system defined by (\ref{eqX}) and (\ref{optimalY}).

\subsection{Proof that $(X_t,Y_t)_{t \geq 0}$ is a coupling}\label{sectionCoupling}

By the previous subsection, we already have the existence of the process $(X_t, Y_t)_{t \geq 0}$ defined as a solution to (\ref{eqX}) and (\ref{optimalY}). However, we still need to show that $(X_t, Y_t)_{t \geq 0}$ is indeed a coupling. If we denote
\begin{equation}\label{defB}
 B(X_{t-},Y_{t-},v,u) := R(X_{t-},Y_{t-})v + (Z_{t-} + v - R(X_{t-},Y_{t-})v)\mathbf{1}_{\{ u < \rho(v, Z_{t-}) \}}
\end{equation}
and
\begin{equation}\label{LtildeDef}
 \widetilde{L}_t := \int_0^t \int_{\{ |v| > m \} \times [0,1]} v N(ds,dv,du) + \int_0^t \int_{\{ |v| \leq m \} \times [0,1]} B(X_{s-},Y_{s-},v,u) \widetilde{N}(ds,dv,du) \,,
\end{equation}
then we can write the equation (\ref{optimalYequation}) for $(Y_t)_{t \geq 0}$ as
\begin{equation*}
 dY_t = b(Y_t)dt + d\widetilde{L}_t \,.
\end{equation*}
Then, if we show that $(\widetilde{L}_t)_{t \geq 0}$ is a L\'{e}vy process with the same finite dimensional distributions as $(L_t)_{t \geq 0}$ defined by (\ref{newLt}), our assertion follows from the uniqueness in law of solutions to the equation (\ref{SDE1}). An analogous fact in the Brownian case was proved using the L\'{e}vy characterization theorem for Brownian motion. Here the proof is more involved, although the idea is very similar. It is sufficient to show two things. First we need to prove that for any $z \in \mathbb{R}^d$ and any $t \geq 0$ we have
\begin{equation}\label{charFunctionsEq}
 \mathbb{E} \exp (i\langle z, \widetilde{L}_t \rangle) = \mathbb{E} \exp (i\langle z, L_t \rangle) \,.
\end{equation}
Then we must also show that for any $t > s \geq 0$ the increment
\begin{equation*}
 \widetilde{L}_t - \widetilde{L}_s
\end{equation*}
is independent of $\mathcal{F}_s$, where $(\mathcal{F}_t)_{t \geq 0}$ is the filtration generated by $(L_t)_{t \geq 0}$.
We will need the following lemma.

\begin{lemma}\label{condLemma}
 Let $f(v,u)$ be a random function on $\{ |v| \leq m \} \times [0,1]$, measurable with respect to $\mathcal{F}_{t_1}$. If
 \begin{equation}\label{condLemmaAssump}
  \mathbb{P} \left(\int_{\{ |v| \leq m \} \times [0,1]} |f(v,u)|^2 \nu(dv) du < \infty \right) = 1 \,,
 \end{equation}
 then
 \begin{equation}\label{condCharFunct}
 \begin{split}
  \mathbb{E} &\left[ \exp \left(i \left\langle z , \int_{t_1}^{t_2} \int_{\{ |v| \leq m \} \times [0,1]} f(v,u) \widetilde{N}(ds,dv,du) \right\rangle \right) \Bigg| \mathcal{F}_{t_1} \right] \\
  &= \exp \left( (t_2-t_1)\int_{\{ |v| \leq m \} \times [0,1]} \left( e^{i \langle z, f(v,u) \rangle} -1 -i\langle z, f(v,u) \rangle \right) \nu(dv) du \right) \,.
  \end{split}
 \end{equation}
\begin{proof}
By a standard argument, if the condition (\ref{condLemmaAssump}) is satisfied, we can approximate $\int_{t_1}^{t_2} \int_{\{ |v| \leq m \} \times [0,1]} f(v,u) \widetilde{N}(ds,dv,du)$ in probability by integrals of step functions $f^n$ of the form
 \begin{equation*}
  f^n(v,u) = \sum_{j=1}^{l_n} c_j \mathbf{1}_{A_j}
 \end{equation*}
 where $A_j$ are pairwise disjoint subsets of $\{ |v| \leq m \} \times [0,1]$ such that $(\nu \times \lambda)(A_j) < \infty$ for all $j$, where $\lambda$ is the Lebesgue measure on $[0,1]$ and $c_j$ are $\mathcal{F}_{t_1}$-measurable random variables. Thus it is sufficient to show (\ref{condCharFunct}) for the step functions $f^n$ and then pass to the limit using the dominated convergence theorem for conditional expectations.
 Indeed, for every $f^n$ we can show that
 \begin{equation*}
  \begin{split}
     \mathbb{E} &\left[ \exp \left( i \left\langle z , \int_{t_1}^{t_2} \int_{\{ |v| \leq m \} \times [0,1]} f^n(v,u) \widetilde{N}(ds,dv,du) \right\rangle \right) \Bigg| \mathcal{F}_{t_1} \right] \\
     &= \mathbb{E} \left[ \prod_{j=1}^{l_n} \exp \left( i \left\langle z , \int_{t_1}^{t_2} \int_{\{ |v| \leq m \} \times [0,1]} c_j \mathbf{1}_{A_j} \widetilde{N}(ds,dv,du) \right\rangle \right) \Bigg| \mathcal{F}_{t_1} \right]\\
     &= \mathbb{E} \left[\prod_{j=1}^{l_n} \exp \left( i \left\langle z , c_j \widetilde{N}((t_1,t_2],A_j) \right\rangle \right) \Bigg| \mathcal{F}_{t_1} \right] \,.
  \end{split}
 \end{equation*}
The random variables $\widetilde{N}((t_1,t_2],A_j)$ are mutually independent and they are all independent of $\mathcal{F}_{t_1}$ and the random variables $c_j$ are $\mathcal{F}_{t_1}$-measurable so we know that we can calculate the above conditional expectation as just an expectation with $c_j$ constant and then plug the random $c_j$ back in. Thus we get
\begin{equation*}
\begin{split}
 \mathbb{E} \prod_{j=1}^{l_n} \exp &\left( i \left\langle z , c_j \widetilde{N}((t_1,t_2],A_j) \right\rangle \right) = \prod_{j=1}^{l_n} \mathbb{E} \exp \left( i \left\langle z , c_j \widetilde{N}((t_1,t_2],A_j) \right\rangle \right)\\
 &=\prod_{j=1}^{l_n} \exp \left( (t_2 - t_1)\left( e^{i\langle z, c_j \rangle}(\nu \times \lambda)(A_j) - 1 - i\langle z, c_j \rangle (\nu \times \lambda)(A_j) \right) \right)\\
 &=\exp \left( (t_2-t_1)\int_{\{ |v| \leq m \} \times [0,1]} \left( e^{i \langle z, f^n(v,u) \rangle} -1 -i\langle z, f^n(v,u) \rangle \right) \nu(dv) du \right) \,,
 \end{split}
\end{equation*}
where in the second step we just used the formula for the characteristic function of the Poisson distribution.
\end{proof}
\end{lemma}

Now we will prove (\ref{charFunctionsEq}) in the special case where
\begin{equation*}
 \widetilde{L}_t = \int_0^t \int_{\{ |v| \leq m \} \times [0,1]} B(X_{s-},Y_{s-},v,u) \widetilde{N}(ds,dv,du) 
\end{equation*}
and the process $(L_t)_{t \geq 0}$ is also considered without the large jumps. Once we have this, it is easy to extend the result to the general case where $(\widetilde{L}_t)_{t \geq 0}$ is given by (\ref{LtildeDef}).

\begin{lemma}\label{charFunctionsTheorem}
 For every $t > 0$ and every $z \in \mathbb{R}^d$ we have
 \begin{equation*}
 \begin{split}
  \mathbb{E}\exp&\left( i\left\langle z, \int_0^t \int_{\{ |v| \leq m \} \times [0,1]} B(X_{s-},Y_{s-},v,u) \widetilde{N}(ds,dv,du)  \right\rangle \right) \\
  &= \mathbb{E}\exp \left( i\left\langle z, \int_0^t \int_{\{ |v| \leq m \} \times [0,1]} v \widetilde{N}(ds,dv,du)  \right\rangle \right) \,.
  \end{split}
 \end{equation*}
\begin{proof}
First recall that we have
\begin{equation*}
 \mathbb{P}\left( \int_{\{ |v| \leq m \} \times [0,1]} |B(X_{t-},Y_{t-},v,u)|^2 \nu(dv) du < \infty \right) = 1
\end{equation*}
(see the remark after the proof of Lemma \ref{linearGrowthLemma}). Then observe that by Lemma \ref{continuityLemma} we know that the square integrated process $B$, i.e., the process
\begin{equation*}
 \int_{\{ |v| \leq m \} \times [0,1]} |B(X_{t-},Y_{t-},v,u)|^2 \nu(dv) du 
\end{equation*}
has left-continuous trajectories. This means that (almost surely) we can approximate $B(X_{t-},Y_{t-},v,u)$ in $L^2([0,t] \times (\{|v| \leq m \}; \nu) \times [0,1])$ by Riemann sums of the form
\begin{equation}\label{RiemannSums}
 B^n(s,v,u) := \sum_{k=0}^{m_n - 1} B(X_{t^n_{k}},Y_{t^n_k},v,u) \mathbf{1}_{(t^n_k,t^n_{k+1}]}(s)
\end{equation}
for some sequence of partitions $0 = t^n_0 < t^n_1 < \ldots < t^n_{m_n} = t$ of the interval $[0,t]$ with the mesh size going to zero as $n \to \infty$. From the general theory of stochastic integration with respect to Poisson random measures (see e.g. \cite{apple}, Section 4.2) it follows that the sequence of integrals $\int_0^t \int_{\{ |v| \leq m \} \times [0,1]} B^n(s,v,u) \widetilde{N}(ds,dv,du)$ converges in probability to the integral $\int_0^t \int_{\{ |v| \leq m \} \times [0,1]} B(X_{s-},Y_{s-},v,u) \widetilde{N}(ds,dv,du)$.
Thus we have
\begin{equation*}
\begin{split}
 \mathbb{E}\exp &\left(i \left\langle z, \int_0^t \int_{\{ |v| \leq m \} \times [0,1]} B^n(s,v,u) \widetilde{N}(ds,dv,du) \right\rangle \right) \\
 &\to \mathbb{E}\exp \left( i \left\langle z, \int_0^t \int_{\{ |v| \leq m \} \times [0,1]} B(X_{s-},Y_{s-},v,u) \widetilde{N}(ds,dv,du) \right\rangle \right)
 \end{split}
\end{equation*}
for any $z \in \mathbb{R}^d$ and $t > 0$, as $n \to \infty$. We will show now that in fact for all $n \in \mathbb{N}$ we have
\begin{equation}\label{equalityOfCharFunctOfRiemannSums}
\begin{split}
 \mathbb{E}\exp &\left( i \left\langle z, \int_0^t \int_{\{ |v| \leq m \} \times [0,1]} B^n(s,v,u) \widetilde{N}(ds,dv,du) \right\rangle \right) \\
 &= \mathbb{E}\exp \left( i \left\langle z,  \int_0^t \int_{\{ |v| \leq m \} \times [0,1]} v \widetilde{N}(ds,dv,du) \right\rangle \right) \,,
 \end{split}
\end{equation}
which will prove the desired assertion. To this end, let us calculate
\begin{equation}\label{conditioningCalculations}
 \begin{split}
  \mathbb{E}&\exp \left(i \left\langle z, \int_0^t \int_{\{ |v| \leq m \} \times [0,1]} B^n(s,v,u) \widetilde{N}(ds,dv,du) \right\rangle \right) \\
  &= \mathbb{E}\exp \left( i \left\langle z, \sum_{k=0}^{m_n - 1}\int_{t^n_k}^{t^n_{k+1}} \int_{\{ |v| \leq m \} \times [0,1]} B(X_{t^n_{k}},Y_{t^n_k},v,u) \widetilde{N}(ds,dv,du) \right\rangle \right)\\
  &=\mathbb{E} \Bigg( \mathbb{E} \Bigg[ \prod_{k=0}^{m_n - 2} \exp \left(i \left\langle z, \int_{t^n_k}^{t^n_{k+1}} \int_{\{ |v| \leq m \} \times [0,1]} B(X_{t^n_{k}},Y_{t^n_k},v,u) \widetilde{N}(ds,dv,du) \right\rangle \right) \\
  &\times \exp \left( i \left\langle z, \int_{t^n_{m_n - 1}}^{t^n_{m_n}} \int_{\{ |v| \leq m \} \times [0,1]} B(X_{t^n_{m_n - 1}},Y_{t^n_{m_n - 1}},v,u) \widetilde{N}(ds,dv,du) \right\rangle \right) \Bigg| \mathcal{F}^n_{m_n - 1}\Bigg] \Bigg) \\
  &=\mathbb{E} \Bigg( \prod_{k=0}^{m_n - 2} \exp \left( i \left\langle z, \int_{t^n_k}^{t^n_{k+1}} \int_{\{ |v| \leq m \} \times [0,1]} B(X_{t^n_{k}},Y_{t^n_k},v,u) \widetilde{N}(ds,dv,du) \right\rangle \right) \\
  &\times \mathbb{E} \left[  \exp \left( i \left\langle z, \int_{t^n_{m_n - 1}}^{t^n_{m_n}} \int_{\{ |v| \leq m \} \times [0,1]} B(X_{t^n_{m_n - 1}},Y_{t^n_{m_n - 1}},v,u) \widetilde{N}(ds,dv,du) \right\rangle \right) \Bigg| \mathcal{F}^n_{m_n - 1}\right]\Bigg)
 \end{split}
\end{equation}
Now we can use Lemma \ref{condLemma} to evaluate the conditional expectation appearing above as
\begin{equation*}
\begin{split}
 &\exp\Bigg( (t^n_{m_n - 1}-t^n_{m_n}) \\
 &\times \int_{\{ |v| \leq m \} \times [0,1]} \left( e^{ i \langle z, B(X_{t^n_{m_n - 1}},Y_{t^n_{m_n - 1}},v,u) \rangle } -1 -i\langle z, B(X_{t^n_{m_n - 1}},Y_{t^n_{m_n - 1}},v,u) \rangle \right) \nu(dv) du \Bigg) \,.
 \end{split}
\end{equation*}
Here comes the crucial part of our proof. We will show that
\begin{equation}\label{crucialPart}
\begin{split}
 \int_{\{ |v| \leq m \} \times [0,1]} &\left( e^{i \langle z, B(X_{t^n_{m_n - 1}},Y_{t^n_{m_n - 1}},v,u) \rangle} -1 -i\langle z, B(X_{t^n_{m_n - 1}},Y_{t^n_{m_n - 1}},v,u) \rangle \right) \nu(dv) du \\
 &= \int_{\{ |v| \leq m \} \times [0,1]} \left( e^{i \langle z, v \rangle} -1 -i\langle z, v \rangle \right) \nu(dv) du \,.
 \end{split}
\end{equation}
Let us fix the values of $X_{t^n_{m_n - 1}}$ and $Y_{t^n_{m_n - 1}}$ for the moment and denote
\begin{equation}\label{RcNotation}
 R:= R(X_{t^n_{m_n - 1}},Y_{t^n_{m_n - 1}}) \text{ and } c:= X_{t^n_{m_n - 1}} - Y_{t^n_{m_n - 1}} = Z_{t^n_{m_n - 1}} \,.
\end{equation}
Then, using the formula (\ref{defB}) we can write
\begin{equation*}
 B(X_{t^n_{m_n - 1}},Y_{t^n_{m_n - 1}},v,u) = Rv + (c + v - Rv)\mathbf{1}_{\{ u < \rho(v, c) \}} \,.
\end{equation*}
Next, integrating over $[0,1]$ with respect to $u$, we get
\begin{equation*}
\begin{split}
 \int_{\{ |v| \leq m \} \times [0,1]} &\left( e^{i \langle z, B(X_{t^n_{m_n - 1}},Y_{t^n_{m_n - 1}},v,u) \rangle} -1 -i\langle z, B(X_{t^n_{m_n - 1}},Y_{t^n_{m_n - 1}},v,u) \rangle \right) \nu(dv) du \\
 &= \int_{\{ |v| \leq m \}} \Bigg( e^{i \langle z, Rv \rangle}\left( e^{i \langle z, c+v-Rv \rangle}\rho(v,c) + (1 - \rho(v,c))\right)-1 \\
 &-i\langle z, Rv \rangle - i\langle z , c + v - Rv \rangle \rho(v,c) \Bigg) \nu(dv) \,.
 \end{split}
\end{equation*}
Since $|B(X_{t^n_{m_n - 1}},Y_{t^n_{m_n - 1}},v,u)|^2$ is integrable with respect to $\nu \times \lambda$ over $\{ |v| \leq m \} \times [0,1]$, 
\begin{equation*}
 e^{i \langle z, Rv \rangle}\left( e^{i \langle z, c+v-Rv \rangle}\rho(v,c) + (1 - \rho(v,c))\right)-1 -i\langle z, Rv \rangle - i\langle z , c + v - Rv \rangle \rho(v,c)
\end{equation*}
is integrable with respect to $\nu$ over $\{ |v| \leq m \}$. Moreover, $e^{i \langle z, Rv \rangle} - 1 -i\langle z, Rv \rangle$ is also integrable over $\{ |v| \leq m \}$. In fact, since $\nu$ is assumed to be rotationally invariant and $R$~is an orthogonal matrix, we easily see that
\begin{equation*}
 \int_{\{ |v| \leq m \}}  \left( e^{i \langle z, Rv \rangle} - 1 -i\langle z, Rv \rangle \right) \nu(dv) = \int_{\{ |v| \leq m \}}  \left( e^{i \langle z, v \rangle} - 1 -i\langle z, v \rangle \right) \nu(dv) \,.
\end{equation*}
We infer that $\left( e^{i \langle z, Rv \rangle} ( e^{i \langle z, c+v-Rv \rangle} - 1 )- i\langle z , c + v - Rv \rangle \right)\rho(v,c)$ is also integrable with respect to $\nu$ over $\{|v| \leq m \}$. Now we will show that the integral of this function actually vanishes. Note that we have $R = I - 2 cc^T / |c|^2$ and since $q$ is the density of a rotationally invariant measure $\nu$, we have $q(Rv) = q(v)$ and $q(Rv - c) = q(v + c)$ for any $v \in \mathbb{R}^d$.
Now
\begin{equation*}
 \begin{split}
  &\int_{\{ |v|\leq m \}} \left( e^{i \langle z, Rv \rangle} ( e^{i \langle z, c+v-Rv \rangle} - 1 )- i\langle z , c + v - Rv \rangle \right)\rho(v,c) \nu(dv) \\
  &= \int_{\{ |v|\leq m \}} \left(  e^{i \langle z, c+v \rangle} - e^{i \langle z, Rv \rangle}- i\langle z , c + v \rangle + i \langle z, Rv \rangle \right)q(v) \wedge q(v+c)\mathbf{1}_{\{|v+c| \leq m \}} dv \\
  &= \int_{\{ |v-c|\leq m \} \cap \{ |v| \leq m \}} \left(  e^{i \langle z, v \rangle} - e^{i \langle z, R(v-c) \rangle}- i\langle z , v \rangle + i \langle z, R(v-c) \rangle \right)q(v-c) \wedge q(v) dv \\
  &= \int_{\{ |v-c| \leq m \} \cap \{ |v| \leq m \}} \left(  e^{i \langle z, v \rangle} - e^{i \langle z, Rv+c \rangle}- i\langle z , v \rangle + i \langle z, Rv+c \rangle \right)q(v-c) \wedge q(v) dv \\
  &= \int_{\{ |Rv-c|\leq m \} \cap \{ |Rv| \leq m \}} \left(  e^{i \langle z, Rv \rangle} - e^{i \langle z, v+c \rangle}- i\langle z , Rv \rangle + i \langle z, v+c \rangle \right)q(Rv-c) \wedge q(Rv) dv \\
  &= \int_{\{ |v+c|\leq m \} \cap \{ |v| \leq m \}} \left(  e^{i \langle z, Rv \rangle} - e^{i \langle z, v+c \rangle}- i\langle z , Rv \rangle + i \langle z, v+c \rangle \right)q(v+c) \wedge q(v) dv \\
  &= \int_{\{ |v| \leq m \}} \left(  e^{i \langle z, Rv \rangle} - e^{i \langle z, v+c \rangle}- i\langle z , Rv \rangle + i \langle z, v+c \rangle \right) \rho(v,c) \nu(dv) \\
  &= - \int_{\{ |v|\leq m \}} \left( e^{i \langle z, Rv \rangle} ( e^{i \langle z, c+v-Rv \rangle} - 1 )- i\langle z , c + v - Rv \rangle \right)\rho(v,c) \nu(dv) \,,
 \end{split}
\end{equation*}
where in the second step we use a change of variables from $v$ to $v-c$, in the third step we use the fact that $Rc = -c$, in the fourth step we change the variables from $v$ to $Rv$ and in the fifth step we use the symmetry properties $|Rv - c| = |v+c|$ and $|Rv| = |v|$. 
Hence we have shown (\ref{crucialPart}). Now we return to our calculations in (\ref{conditioningCalculations}) and compute
\begin{equation}\label{endOfFirstStep}
 \begin{split}
  &\mathbb{E} \Bigg( \prod_{k=0}^{m_n - 2} \exp \left( i \left\langle z, \int_{t^n_k}^{t^n_{k+1}} \int_{\{ |v| \leq m \} \times [0,1]} B(X_{t^n_{k}},Y_{t^n_k},v,u) \widetilde{N}(ds,dv,du) \right\rangle \right) \\
  &\times \mathbb{E} \left[  \exp \left( i \left\langle z, \int_{t^n_{m_n - 1}}^{t^n_{m_n}} \int_{\{ |v| \leq m \} \times [0,1]} B(X_{t^n_{m_n - 1}},Y_{t^n_{m_n - 1}},v,u) \widetilde{N}(ds,dv,du) \right\rangle \right) \Bigg| \mathcal{F}^n_{m_n - 1}\right]\Bigg)\\
  &= \exp\left( (t^n_{m_n}-t^n_{m_n-1})\int_{\{ |v| \leq m \} \times [0,1]} \left( e^{i \langle z, v \rangle} -1 -i\langle z, v \rangle \right) \nu(dv) du \right) \\
  &\times \mathbb{E} \left( \prod_{k=0}^{m_n - 2} \exp \left( i \left\langle z, \int_{t^n_k}^{t^n_{k+1}} \int_{\{ |v| \leq m \} \times [0,1]} B(X_{t^n_{k}},Y_{t^n_k},v,u) \widetilde{N}(ds,dv,du) \right\rangle \right) \right) \,.
 \end{split}
\end{equation}
Then we can just repeat all the steps from (\ref{conditioningCalculations}) to (\ref{endOfFirstStep}), this time conditioning on $\mathcal{F}^n_{m_n-2}$, 
and after repeating this procedure $m_n - 1$ times, we get (\ref{equalityOfCharFunctOfRiemannSums}).
\end{proof}
\end{lemma}

It remains now to show the independence of the increments of $(\widetilde{L}_t)_{t \geq 0}$.

\begin{lemma}\label{independenceLemma}
 Under the above assumptions, for any $t_2 > t_1 \geq 0$ the random variable $\widetilde{L}_{t_2} -\widetilde{L}_{t_1}$ is independent of $\mathcal{F}_{t_1}$.
 \begin{proof}
   We will show that for an arbitrary $\mathcal{F}_{t_1}$-measurable random variable $\xi$ and for any $z_1$,~$z_2 \in~\mathbb{R}^d$ we have
  \begin{equation*}
  \begin{split}
   &\mathbb{E}\exp\left(i\left\langle z_1, \int_{t_1}^{t_2} \int_{\{ |v| \leq m \} \times [0,1]} B(X_{s-},Y_{s-},v,u) \widetilde{N}(ds,dv,du) \right\rangle + i\langle z_2 , \xi \rangle \right)\\
   &= \mathbb{E}\exp\left(i\left\langle z_1, \int_{t_1}^{t_2} \int_{\{ |v| \leq m \} \times [0,1]} B(X_{s-},Y_{s-},v,u) \widetilde{N}(ds,dv,du) \right\rangle \right) \cdot \mathbb{E} \exp ( i\langle z_2 , \xi \rangle ) \,.
    \end{split}
  \end{equation*}
As in the proof of Lemma \ref{charFunctionsTheorem}, the integral $\int_{t_1}^{t_2} \int_{\{ |v| \leq m \} \times [0,1]} B(X_{s-},Y_{s-},v,u) \widetilde{N}(ds,dv,du)$ can be approximated by integrals of Riemann sums $B^n(s,v,u)$ that have been defined by (\ref{RiemannSums}) for some sequence of partitions
$t_1 = t^n_0 < t^n_1 < \ldots < t^n_{m_n} = t_2$ such that $\delta_n := \max_{k \in \{ 0, \ldots , m_n - 1\}} |t^n_{k+1} - t^n_k| \to 0$ as $n \to \infty$. Denote
\begin{equation*}
 I^n_k := \int_{t^n_k}^{t^n_{k+1}} \int_{\{ |v| \leq m \} \times [0,1]} B(X_{t^n_k},Y_{t^n_k},v,u) \widetilde{N}(ds,dv,du) \,, \text{  } I^n := \sum_{k=0}^{m_n - 1} I^n_k \,.
\end{equation*}
Then we have
\begin{equation}\label{lemmaIndependenceCalculations}
\begin{split}
 \mathbb{E} \exp &\left( i \langle z_1 , I^n \rangle + i\langle z_2, \xi \rangle \right) = \mathbb{E} \left( \exp (i \langle z_2, \xi \rangle) \prod_{k=0}^{m_n - 1} \exp (i \langle z_1 , I^n_k \rangle) \right) \\
 &= \mathbb{E} \left( \mathbb{E} \left[ \exp (i \langle z_2, \xi \rangle) \prod_{k=0}^{m_n - 1} \exp (i \langle z_1 , I^n_k \rangle) \Bigg| \mathcal{F}_{t^n_{m_n - 1}} \right] \right) \\
 &= \mathbb{E} \left( \exp (i \langle z_2, \xi \rangle) \prod_{k=0}^{m_n - 2} \exp (i \langle z_1 , I^n_k \rangle) \mathbb{E} \left[ \exp (i \langle z_1 , I^n_{m_n - 1} \rangle) \Bigg| \mathcal{F}_{t^n_{m_n - 1}} \right] \right) \,,
 \end{split}
\end{equation}
where in the last step we used the fact that for every $k \in \{ 0, \ldots, m_n - 1 \}$ the random variable $\xi$ is $\mathcal{F}_{t_1} \subset \mathcal{F}_{t^n_k}$-measurable. Now, using Lemma \ref{condLemma} and our calculations from the proof of Lemma \ref{charFunctionsTheorem}, we can show that
\begin{equation*}
 \mathbb{E} \left[ \exp (i \langle z_1 , I^n_{m_n - 1} \rangle) \Bigg| \mathcal{F}_{t^n_{m_n - 1}} \right] = \mathbb{E} \exp \left( i \left\langle z_1 , \int_{t^n_{m_n - 1}}^{t^n_{m_n}} \int_{\{ |v| \leq m \} \times [0,1]} v \widetilde{N}(ds,dv,du) \right\rangle \right)
\end{equation*}
and thus we see that the expression on the right hand side of (\ref{lemmaIndependenceCalculations}) is equal to
\begin{equation*}
\mathbb{E} \exp \left( i \left\langle z_1 , \int_{t^n_{m_n - 1}}^{t^n_{m_n}} \int_{\{ |v| \leq m \} \times [0,1]} v \widetilde{N}(ds,dv,du) \right\rangle \right) \mathbb{E} \left( \exp (i \langle z_2, \xi \rangle) \prod_{k=0}^{m_n - 2} \exp (i \langle z_1 , I^n_k \rangle) \right) \,.
\end{equation*}
Thus, by repeating the above procedure $m_n - 1$ times (conditioning on the consecutive $\sigma$-fields $\mathcal{F}_{t^n_k}$), we get
\begin{equation}\label{eqIndependenceProof}
\begin{split}
  \mathbb{E} \exp \left( i \langle z_1 , I^n \rangle + i\langle z_2, \xi \rangle \right) &= \mathbb{E} \exp (i \langle z_2 , \xi \rangle) \\
  &\times \prod_{k=0}^{m_n - 1} \mathbb{E} \exp \left( i \left\langle z_1 , \int_{t^n_{k}}^{t^n_{k+1}} \int_{\{ |v| \leq m \} \times [0,1]} v \widetilde{N}(ds,dv,du) \right\rangle \right) \,.
  \end{split}
\end{equation}
However, by the same argument as above we can show that
\begin{equation*}
 \prod_{k=0}^{m_n - 1} \mathbb{E} \exp \left( i \left\langle z_1 , \int_{t^n_{k}}^{t^n_{k+1}} \int_{\{ |v| \leq m \} \times [0,1]} v \widetilde{N}(ds,dv,du) \right\rangle \right) = \mathbb{E} \exp ( i \langle z_1, I^n \rangle) \,.
\end{equation*}
Since $I^n$ converges in probability to $\int_{t_1}^{t_2} \int_{\{ |v| \leq m \} \times [0,1]} B(X_{s-},Y_{s-},v,u) \widetilde{N}(ds,dv,du)$, we get
\begin{equation*}
 \mathbb{E} \exp (i \langle z_1 , I^n \rangle) \to \mathbb{E} \exp \left( i \left\langle z_1 , \int_{t_1}^{t_2} \int_{\{ |v| \leq m \} \times [0,1]} B(X_{s-},Y_{s-},v,u) \widetilde{N}(ds,dv,du) \right\rangle \right)
\end{equation*}
and, by passing to a subsequence for which almost sure convergence holds and using the dominated convergence theorem, we get
\begin{equation*}
\begin{split}
 \mathbb{E} \exp &(i \langle z_1 , I^n \rangle + i\langle z_2 , \xi \rangle) \\
 &\to \mathbb{E} \exp \left( i \left\langle z_1 , \int_{t_1}^{t_2} \int_{\{ |v| \leq m \} \times [0,1]} B(X_{s-},Y_{s-},v,u) \widetilde{N}(ds,dv,du) \right\rangle + i \langle z_2, \xi \rangle \right) \,,
 \end{split}
\end{equation*}
which proves the desired assertion.
 \end{proof}
\end{lemma}

\section{Proof of the inequality (\ref{mainInequality})}\label{Section3}

In this section we want to apply the coupling that we constructed in Section \ref{Section2} to prove Corollary \ref{mainTheorem}, which follows easily from the inequality (\ref{mainInequality}). Namely, in order to obtain
\begin{equation}\label{WfInequality}
 W_f(\mu p_t, \nu p_t) \leq e^{-ct}W_f(\mu, \nu) \,,
\end{equation}
we will prove that
\begin{equation}\label{expectationInequality}
 \mathbb{E}f(|X_t - Y_t|) \leq e^{-ct}\mathbb{E}f(|X_0 - Y_0|) \,,
\end{equation}
where $(X_t,Y_t)_{t \geq 0}$ is the coupling defined by (\ref{eqX}) and (\ref{optimalY}) and the laws of the random variables $X_0$ and $Y_0$ are $\mu$ and $\nu$, respectively. Obviously, straight from the definition of the distance $W_f$ we see that for any coupling $(X_t,Y_t)_{t \geq 0}$ the expression $\mathbb{E}f(|X_t - Y_t|)$ gives an upper bound for $W_f(\mu p_t, \nu p_t)$ and since we can prove (\ref{expectationInequality}) for any coupling of the initial conditions $X_0$ and $Y_0$, it is easy to see that (\ref{expectationInequality}) indeed implies (\ref{WfInequality}). Note that without loss of generality we can assume that $\mathbb{P}(X_0 \neq Y_0) = 1$. Indeed, given any probability measures $\mu$ and $\nu$ we can decompose them by writing
\begin{equation}\label{measureDecomposition}
 \mu = \mu \wedge \nu + \widetilde{\mu} \text{  and } \nu = \mu \wedge \nu + \widetilde{\nu}
\end{equation}
for some finite measures $\widetilde{\mu}$ and $\widetilde{\nu}$ on $\mathbb{R}^d$. Then, if $\alpha := (\mu \wedge \nu)(\mathbb{R}^d) \in (0,1)$, we can define probability measures $\overline{\mu} := \widetilde{\mu} / \widetilde{\mu}(\mathbb{R}^d)$ and $\overline{\nu} := \widetilde{\nu} / \widetilde{\nu}(\mathbb{R}^d)$ and we can easily show that $W_f(\mu, \nu) = (1-\alpha)W_f(\overline{\mu}, \overline{\nu})$. Obviously, the decomposition (\ref{measureDecomposition}) is preserved by the semigroup $(p_t)_{t \geq 0}$ and thus we see that in order to show (\ref{WfInequality}) it is sufficient to show that $W_f(\overline{\mu}p_t, \overline{\nu}p_t) \leq e^{-ct} W_f(\overline{\mu}, \overline{\nu})$. 

In our proof we will aim to obtain estimates of the form
\begin{equation}\label{preGronwall}
 \mathbb{E}f(|Z_t|) - \mathbb{E}f(|Z_0|) \leq \mathbb{E} \int_0^t -cf(|Z_s|)ds \,,  
\end{equation}
for some constant $c > 0$, where $Z_t = X_t - Y_t$, which by the Gronwall inequality will give us (\ref{expectationInequality}). We assume that $f$ is of the form
\begin{equation*}
 f = f_1 + f_2 \,,
\end{equation*}
where $f_1 \in \mathcal{C}^2$, $f_1' \geq 0$, $f_1'' \leq 0$ and $f_1(0)=0$ and $f_2 = a\mathbf{1}_{(0,\infty)}$ for some constant $a > 0$ to be chosen later. We also choose $f_1$ in such a way that $f_1'(0) = 1$ and thus $f_1' \leq 1$ since $f_1'$~is decreasing. Recall that our coupling is defined in such a way that the equation for the difference process $Z_t = X_t - Y_t$ is given by
\begin{equation}\label{optimalZ}
 \begin{split}
  dZ_t &= (b(X_t) - b(Y_t))dt + \int_{\{ |v| \leq m \} \times [0,1]} (I-R(X_{t-},Y_{t-}))v \widetilde{N}(dt,dv,du) \\
  &- \int_{\{ |v| \leq m \} \times [0,1]} (Z_{t-} + v - R(X_{t-},Y_{t-})v)\mathbf{1}_{\{ u < \rho(v, Z_{t-}) \}} \widetilde{N}(dt,dv,du) \,.
 \end{split}
\end{equation}
Note that the jumps of size greater than $m$ cancel out, since we apply synchronous coupling for $|v| > m$ in our construction of the process $(Y_t)_{t \geq 0}$. In order to simplify the notation, let us denote
\begin{equation}\label{defBigA}
 A(X_{t-},Y_{t-},v,u):= -(Z_{t-} + v - R(X_{t-},Y_{t-})v)\mathbf{1}_{\{ u < \rho(v, Z_{t-}) \}} \,.
\end{equation}
Then we can write
\begin{equation}\label{optimalZwithA}
 \begin{split}
  dZ_t &= (b(X_t) - b(Y_t))dt + \int_{\{ |v| \leq m \} \times [0,1]} (I-R(X_{t-},Y_{t-}))v \widetilde{N}(dt,dv,du) \\
  &+ \int_{\{ |v| \leq m \} \times [0,1]} A(X_{t-},Y_{t-},v,u) \widetilde{N}(dt,dv,du) \,.
 \end{split}
\end{equation}
Let us split our computations into two parts by writing
\begin{equation}\label{ItoPreliminary}
  \mathbb{E}f(|Z_t|) -  \mathbb{E}f(|Z_0|) = \mathbb{E}f_1(|Z_t|)  -  \mathbb{E}f_1(|Z_0|) +  a\mathbb{E}\mathbf{1}_{(0,\infty)}(|Z_t|) -  a\mathbb{E}\mathbf{1}_{(0,\infty)}(|Z_0|) \,.
\end{equation}
We will first deal with finding an appropriate formula for $f_1$ by bounding the difference $\mathbb{E}f_1(|Z_t|)  -  \mathbb{E}f_1(|Z_0|)$ from above. This way we will obtain some estimates that are valid only under the assumption that $|Z_s| > \delta$ for some $\delta > 0$ and all $s \in [0,t]$. We will then use the discontinuous part $f_2$ of our distance function $f$ to improve these results and obtain bounds that hold regardless of the value of $|Z_s|$. We will start the proof by applying the It\^{o} formula for L\'{e}vy processes (see e.g. \cite{apple}, Theorem 4.4.10) to the equation (\ref{optimalZwithA}) and the function $g(x):=f_1(|x|)$. We have
\begin{equation}\label{secondPartialg}
 \partial_ig(x) = f_1'(|x|) \frac{x_i}{|x|} \text{   and  } \partial_j \partial_i g(x) = f_1''(|x|)\frac{x_j x_i}{|x|^2} + f_1'(|x|)\left(\delta_{ij}\frac{1}{|x|} - \frac{x_j x_i}{|x|^3}\right) \,,
\end{equation}
where $\delta_{ij}$ is the Kronecker delta. By the It\^{o} formula we have
\begin{equation}\label{ItoForg}
 g(Z_t) - g(Z_0) = \sum_{i=1}^{d} \int_0^t \partial_ig(Z_{s-})dZ^i_s + \sum_{s \in (0,t]} \left( g(Z_s) - g(Z_{s-}) - \sum_{i=1}^{d} \partial_ig(Z_{s-}) \Delta Z_s^i \right) \,,
\end{equation}
where $Z_t = (Z_t^1, \ldots , Z_t^d)$ and $\Delta Z_t = Z_t - Z_{t-}$. Using the Taylor formula we can write
\begin{equation*}
  g(Z_s) - g(Z_{s-}) - \sum_{i=1}^{d} \partial_ig(Z_{s-}) \Delta Z_s^i = \sum_{i,j = 1}^d \int_0^1 (1-u) \partial_j \partial_i g(Z_{s-} + u \Delta Z_s)du \Delta Z_s^i \Delta Z_s^j \,.
\end{equation*}
Denoting $W_{s,u}:=Z_{s-} + u \Delta Z_s$ and using (\ref{secondPartialg}), we can further evaluate the above expression as
\begin{equation}\label{TaylorSum}
 \sum_{i,j = 1}^d \int_0^1 (1-u) \left[ f_1''(|W_{s,u}|)\frac{W_{s,u}^j W_{s,u}^i}{|W_{s,u}|^2} + f_1'(|W_{s,u}|)\frac{1}{|W_{s,u}|}\left( \delta_{ij} - \frac{W_{s,u}^j W_{s,u}^i}{|W_{s,u}|^2} \right) \right] du \Delta Z_s^i \Delta Z_s^j \,.
\end{equation}
Observe now that for every $s \in (0,t]$ and every $u \in (0,1)$ the vectors $\Delta Z_s$ and $W_{s,u}$ are parallel. This follows from the fact that if $\Delta Z_s \neq 0$ (i.e., there is a jump at $s$) then $Y_s$ is equal either to $X_s$ or to $R(X_{s-},Y_{s-})X_s$ and hence $Z_s$ is equal either to zero or to $2e_{s-}e_{s-}^TX_s$, which is obviously parallel to $Z_{s-}$. Thus we always have
\begin{equation*}
 \sum_{i=1}^d W_{s,u}^i \Delta Z_s^i = \langle W_{s,u} , \Delta Z_s \rangle = \pm |W_{s,u}| \cdot |\Delta Z_s|
\end{equation*}
and in consequence (\ref{TaylorSum}) is equal to
\begin{equation*}
\begin{split}
 \int_0^1 &(1-u) \left[ f_1''(|W_{s,u}|)\frac{|W_{s,u}|^2 |\Delta Z_s|^2}{|W_{s,u}|^2} + f_1'(|W_{s,u}|)\frac{1}{|W_{s,u}|}\left( |\Delta Z_s|^2 - \frac{|W_{s,u}|^2 |\Delta Z_s|^2}{|W_{s,u}|^2} \right) \right] du \\
 &= \int_0^1 (1-u) f_1''(|W_{s,u}|) |\Delta Z_s|^2 du \,,
 \end{split}
\end{equation*}
so we see that the second sum in (\ref{ItoForg}) is of the form
\begin{equation*}
 \sum_{s \in (0,t]}  \left( |\Delta Z_s|^2 \int_0^1 (1-u) f_1''(|Z_{s-} + u \Delta Z_s|) du \right) \,.
\end{equation*}
Hence we can write (\ref{ItoForg}) as
\begin{equation}\label{ItoOptimal}
 \begin{split}
 f_1(|Z_t|) - f_1(|Z_0|) &= \int_0^t f_1'(|Z_{s-}|)\frac{1}{|Z_{s-}|}\langle Z_{s-} , b(X_{s-}) - b(Y_{s-}) \rangle ds \\
 &+\int_0^t \int_{\{ |v| \leq m \} \times [0,1]} f_1'(|Z_{s-}|)\frac{1}{|Z_{s-}|}\langle Z_{s-},(I-R(X_{s-},Y_{s-}))v \rangle \widetilde{N}(ds,dv,du)  \\
 &+\int_0^t \int_{\{ |v| \leq m \} \times [0,1]} f_1'(|Z_{s-}|)\frac{1}{|Z_{s-}|}\langle Z_{s-},A(X_{s-},Y_{s-},v,u) \rangle \widetilde{N}(ds,dv,du)\\
 &+\sum_{s \in (0,t]}  \left( |\Delta Z_s|^2 \int_0^1 (1-u) f_1''(|Z_{s-} + u \Delta Z_s|) du \right) \,.
 \end{split}
\end{equation}
Note that the above formula holds only for $t < T$, where $T$ is the coupling time defined by (\ref{theCouplingTime}). However, for $t \geq T$ we have $Z_t = 0$ so if we want to obtain (\ref{expectationInequality}), it is sufficient to bound $\mathbb{E}f(|Z_{t \wedge T}|)$. In order to calculate the expectations of the above terms we will use a sequence of stopping times $(\tau_n)_{n=1}^{\infty}$ defined by
\begin{equation*}
 \tau_n := \inf \{ t \geq 0 : |Z_t| \notin (1/n , n) \} \,.
\end{equation*}
Note that we have $\tau_n \to T$ as $n \to \infty$, which follows from non-explosiveness of $(Z_t)_{t \geq 0}$, which in turn is a consequence of non-explosiveness of the solution to (\ref{SDE1}). Now we will split our computations into several lemmas.

\begin{lemma}\label{lemma31}
We have
\begin{equation*}
 \mathbb{E} \int_0^{t \wedge \tau_n} \int_{\{ |v| \leq m \} \times [0,1]} f_1'(|Z_{s-}|)\frac{1}{|Z_{s-}|}\langle Z_{s-},(I-R(X_{s-},Y_{s-}))v \rangle \widetilde{N}(ds,dv,du) = 0 \,.
 \end{equation*}
 \begin{proof}
  Observe that
\begin{equation*}
 \langle Z_{s-} , (I-R(X_{s-},Y_{s-}))v \rangle = \langle Z_{s-} , 2e_{s-}e_{s-}^T v \rangle = 2 \langle e_{s-}, v \rangle \langle Z_{s-} , \frac{Z_{s-}}{|Z_{s-}|} \rangle = 2\langle e_{s-}, v \rangle |Z_{s-}|
\end{equation*}
and therefore
\begin{equation*}
\begin{split}
 \int_0^{t \wedge \tau_n} \int_{\{ |v| \leq m \} \times [0,1]} f_1'(|Z_{s-}|) \frac{1}{|Z_{s-}|}&\langle Z_{s-} , (I-R(X_{s-},Y_{s-}))v \rangle \widetilde{N}(ds,dv,du) \\
 &= 2 \int_0^{t \wedge \tau_n} \int_{\{ |v| \leq m \} \times [0,1]} f_1'(|Z_{s-}|)\langle e_{s-}, v \rangle \widetilde{N}(ds,dv,du) \,.
 \end{split}
\end{equation*}
By the Cauchy-Schwarz inequality and the fact that $f_1' \leq 1$, for any $t \geq 0$ we have
\begin{equation*}
 \int_0^{t \wedge \tau_n} \int_{\{ |v| \leq m\} \times [0,1]} |f_1'(|Z_{s-}|)|^2|\langle e_{s-}, v \rangle|^2 \nu(dv) du ds \leq \int_0^{t} \int_{\{ |v| \leq m\} \times [0,1]} |v|^2 \nu(dv) du ds < \infty \,,
\end{equation*}
which implies that
\begin{equation*}
 \int_0^{t \wedge \tau_n} \int_{\{ |v| \leq m \} \times [0,1]} f_1'(|Z_{s-}|)\langle e_{s-} , v \rangle \widetilde{N}(ds,dv,du)
\end{equation*}
is a martingale, from which we immediately obtain our assertion.
 \end{proof}
\end{lemma}

\begin{lemma}\label{lemma32}
We have
 \begin{equation*}
  \mathbb{E} \int_0^{t \wedge \tau_n} \int_{\{ |v| \leq m \} \times [0,1]} f_1'(|Z_{s-}|)\frac{1}{|Z_{s-}|} \langle Z_{s-},A(X_{s-},Y_{s-},v,u) \rangle \widetilde{N}(ds,dv,du) = 0 \,.
 \end{equation*}
\begin{proof}
 By the Cauchy-Schwarz inequality and the fact that $f_1' \leq 1$, we have
\begin{equation*}
\begin{split}
  &\int_0^{t \wedge \tau_n} \int_{\{ |v| \leq m \} \times [0,1]} |f_1'(|Z_{s-}|)\frac{1}{|Z_{s-}|}\langle Z_{s-},A(X_{s-},Y_{s-},v,u) \rangle|^2 \nu(dv) du ds \\
 &\leq   \int_0^{t \wedge \tau_n} \int_{\{ |v| \leq m \} \times [0,1]} |A(X_{s-},Y_{s-},v,u)|^2 \nu(dv) du ds \,.
 \end{split}
\end{equation*}
Using the bounds obtained in Lemma \ref{linearGrowthLemma} and the fact that $|Z_s| \leq n$ for $s \leq \tau_n$, we can bound the integral above by a constant. Thus we see that the process
\begin{equation*}
 \int_0^{t \wedge \tau_n} \int_{\{ |v| \leq m \} \times [0,1]} f_1'(|Z_{s-}|)\frac{1}{|Z_{s-}|}\langle Z_{s-},A(X_{s-},Y_{s-},v,u) \rangle \widetilde{N}(ds,dv,du) 
\end{equation*}
is a martingale.
\end{proof}

\end{lemma}

\begin{lemma}\label{lemma33}
 For any $t > 0$, we have
 \begin{equation*}
  \mathbb{E} \sum_{s \in (0,t]}  \left( |\Delta Z_s|^2 \int_0^1 (1-u) f_1''(|Z_{s-} + u \Delta Z_s|) du \right) \leq C_{\varepsilon} \mathbb{E} \int_0^t \bar{f}_{\varepsilon}(|Z_{s-}|) \mathbf{1}_{\{ |Z_{s-}| > \delta \}} ds \,,
 \end{equation*}
where $\delta > 0$, $\varepsilon \leq \delta \wedge 2m$, the constant $C_{\varepsilon}$ is defined by
\begin{equation*}
 C_{\varepsilon} := 2 \int_{ -\varepsilon/4  }^0 |y|^2 \nu_1(dy) \,,
\end{equation*}
where $\nu_1$ is the first marginal of $\nu$ and the function $\bar{f}_{\varepsilon}$ is defined by
\begin{equation*}
 \bar{f}_{\varepsilon}(y):= \sup_{x \in (y - \varepsilon, y)}f_1''(x) \,.
\end{equation*}

\begin{remark}\label{remarkChoiceEpsilonDelta}
 Note that the above estimate holds for any $\delta > 0$ and $\varepsilon \leq \delta \wedge 2m$ as long as $\varepsilon$ satisfies Assumption \ref{Assumption4} and $m$ is sufficiently large (see (\ref{choiceOfm1}) below). Even though our calculations from the proof of Lemma \ref{lemmaDiscontinuous} indicate that later on we should choose $\delta$ and $m$ to be the constants from Assumption \ref{Assumption3}, here in Lemma \ref{lemma33} we do not use the condition (\ref{overlapCondition}). Note that if the condition (\ref{overlapCondition}) from Assumption \ref{Assumption3} is satisfied by more than one value of $\delta$ (which is the case for most typical examples), there appears a question of the optimal choice of $\delta$ and $\varepsilon$ that would maximize the contractivity constant $c$ defined by (\ref{contractivityConstant}) via (\ref{defC}) and (\ref{choiceOfm2}). The answer to this depends on the particular choice of the noise $(L_t)_{t \geq 0}$. It is non-trivial though, even in simple cases, since $c$ depends on $\delta$ and $\varepsilon$ in a convoluted way (see the discussion in Example \ref{exampleConstantBounds}).
\end{remark}

\begin{remark}\label{remarkInterval}
 In the proof of the inequality (\ref{mainInequality}), if we want to obtain an inequality of the form (\ref{preGronwall}) from (\ref{ItoOptimal}), we need to bound the sum appearing in (\ref{ItoOptimal}) by a strictly negative term. For technical reasons that will become apparent in the proof of Lemma \ref{lemmaFunctionalInequality} (see the remarks after (\ref{hereWeUseDecreasingDerivative})), we will use the supremum of the second derivative of $f_1$ over ``small'' jumps that decrease the distance between $X_t$ and $Y_t$.
\end{remark}

\begin{proof}
 Observe that for every $u \in (0,1)$ we have
 \begin{equation}\label{fbisBoundBySup}
\begin{split}
 f_1''(|Z_{s-} + u \Delta Z_s|) &= f_1''(|Z_{s-} + u \Delta Z_s|)(\mathbf{1}_{\{ |Z_s| \in (|Z_{s-}|-\varepsilon,|Z_{s-}|) \}} + \mathbf{1}_{\{  |Z_s| \notin (|Z_{s-}|-\varepsilon,|Z_{s-}|)  \} }) \\
 &\leq \sup_{x \in (|Z_{s-}| - \varepsilon, |Z_{s-}|)}f_1''(x)\mathbf{1}_{\{ |Z_s| \in (|Z_{s-}|-\varepsilon,|Z_{s-}|) \}} \,.
 \end{split}
\end{equation}
Indeed, $f_1$ is assumed to be concave, and thus $f_1''$ is negative, so
\begin{equation*}
 f_1''(|Z_{s-} + u \Delta Z_s|)\mathbf{1}_{\{  |Z_s| \notin (|Z_{s-}|-\varepsilon,|Z_{s-}|)  \} } \leq 0 \,.
\end{equation*}
We also know that the vectors $Z_{s-}$ and $\Delta Z_s$ are parallel, hence if $|Z_s| \in (|Z_{s-}|-\varepsilon,|Z_{s-}|)$, then $|Z_{s-} + u \Delta Z_s| = |Z_{s-}| - u|\Delta Z_s|$ for all $u \in (0,1)$. In particular, we have $|\Delta Z_s| \in (0, \varepsilon)$ and $|Z_{s-} + u \Delta Z_s| \in (|Z_{s-}|-\varepsilon,|Z_{s-}|)$ for all $u \in (0,1)$ and hence we have (\ref{fbisBoundBySup}).

Now let $\delta > 0$ be a positive constant (as mentioned in Remark \ref{remarkChoiceEpsilonDelta}, it can be the constant from Assumption \ref{Assumption3}). Here we introduce an additional factor involving $\delta$ in order for the integral in (\ref{integralCepsilonDelta}) to be bounded from below by a positive constant. We have
\begin{equation*}
 \sup_{x \in (y - \varepsilon, y)}f_1''(x) \cdot \mathbf{1}_{\{ |y| > \delta \}} \geq \sup_{x \in (y - \varepsilon, y)}f_1''(x) \,,
\end{equation*}
so we can write
\begin{equation}\label{JumpsSumEstimate}
\begin{split}
 \sum_{s \in (0,t]} & \left( |\Delta Z_s|^2 \int_0^1 (1-u) f_1''(|Z_{s-} + u \Delta Z_s|) du \right) \\
 &\leq  \sum_{s \in (0,t]}  \left( \frac{1}{2}|\Delta Z_s|^2  \bar{f}_{\varepsilon}(|Z_{s-}|)  \right)\mathbf{1}_{\{ |Z_s| \in (|Z_{s-}|-\varepsilon,|Z_{s-}|) \}} \mathbf{1}_{\{ |Z_{s-}| > \delta \}} \,.
\end{split}
\end{equation}
Now observe that
\begin{equation*}
 \{ |Z_s| \in (|Z_{s-}|-\varepsilon,|Z_{s-}|) \} = \{ |Z_s| < |Z_{s-}|\} \cap \{ |\Delta Z_s| < \varepsilon \} \,,
\end{equation*}
and the condition $|Z_s| < |Z_{s-}|$ is equivalent to $\langle \Delta Z_s , 2Z_{s-}+ \Delta Z_s \rangle < 0$, so we have
\begin{equation*}
 \mathbf{1}_{\{ |Z_s| \in (|Z_{s-}|-\varepsilon,|Z_{s-}|) \}} = \mathbf{1}_{\{ |\Delta Z_s| < \varepsilon \}}\mathbf{1}_{\{ \langle \Delta Z_s , 2Z_{s-}+ \Delta Z_s \rangle < 0 \}} \,.
\end{equation*}
Now we can use the equation (\ref{optimalZ}) describing the dynamics of the jumps of the process $(Z_t)_{t \geq 0}$ and express the sum on the right hand side of (\ref{JumpsSumEstimate}) as an integral with respect to the Poisson random measure $N$ associated with $(L_t)_{t \geq 0}$. However, since all the terms in this sum are negative, we can additionally bound it from above by a sum taking into account only the jumps for which $u \geq \rho(v,Z_{s-})$, i.e., only the reflected jumps. After doing all this, we get
\begin{equation*}
\begin{split}
  \mathbb{E} \sum_{s \in (0,t]}  &\left( \frac{1}{2}|\Delta Z_s|^2  \bar{f}_{\varepsilon}(|Z_{s-}|)  \right)\mathbf{1}_{\{ |Z_s| \in (|Z_{s-}|-\varepsilon,|Z_{s-}|) \}}  \mathbf{1}_{\{ |Z_{s-}| > \delta \}} \\
  &\leq \frac{1}{2} \mathbb{E} \int_0^t \int_{\{ |v| \leq m \} \times [0,1]} |2e_{s-}e_{s-}^Tv|^2 \bar{f}_{\varepsilon}(|Z_{s-}|) \mathbf{1}_{\{ |2e_{s-}e_{s-}^Tv|<\varepsilon \}} \\ &\times \mathbf{1}_{\{ \langle 2e_{s-}e_{s-}^Tv , 2Z_{s-} + 2e_{s-}e_{s-}^Tv \rangle < 0 \} } \mathbf{1}_{\{ |Z_{s-}| > \delta \}} N(ds,dv,du) \,.
\end{split}
\end{equation*}
Note that
\begin{equation*}
\begin{split}
 \langle e_{s-}e_{s-}^Tv, Z_{s-} + e_{s-}e_{s-}^Tv \rangle &= \left\langle \langle e_{s-},v\rangle e_{s-}, |Z_{s-}|e_{s-} + \langle e_{s-},v\rangle e_{s-} \right\rangle \\
 &= \langle e_{s-},v \rangle \left( |Z_{s-}| + \langle e_{s-},v \rangle \right)
\end{split}
\end{equation*}
and thus we can express the expectation above as
\begin{equation*}
\begin{split}
 2\mathbb{E} &\int_0^t \bar{f}_{\varepsilon}(|Z_{s-}|) \int_{\{ |v| \leq m \} \times [0,1]} |\langle e_{s-},v \rangle|^2  \mathbf{1}_{\{ |\langle e_{s-},v \rangle| < \varepsilon/2 \}} \\
 &\times \mathbf{1}_{\{ \langle e_{s-},v \rangle \left( |Z_{s-}| + \langle e_{s-},v \rangle \right)  < 0 \} } \mathbf{1}_{\{ |Z_{s-}| > \delta \}} \nu(dv) du ds \,.
\end{split}
\end{equation*}
Now denote $\nu^m(dv) := \mathbf{1}_{\{ |v| \leq m \}} \nu(dv)$ and observe that if we consider the image $\nu^m \circ h_w^{-1}$ of the measure $\nu^m$ by the mapping $h_w: \mathbb{R}^d \to \mathbb{R}$ defined by $h_w(v) = \langle w, v \rangle$ for a unit vector $w \in \mathbb{R}^d$, then due to the rotational invariance of $\nu^m$, the measure $\nu^m \circ h_v^{-1}$ is independent of the choice of $w$, i.e.,
\begin{equation*}
 \nu^m \circ h_w^{-1} = \nu^m_1 \text{ for all unit vectors } w \in \mathbb{R}^d \,,  
\end{equation*}
where $\nu^m_1$ is the first marginal of $\nu^m$ (and therefore it is the jump measure of a one-dimensional L\'{e}vy process being a projection of $(L_t)_{t \geq 0}$ with truncated jumps, see e.g. \cite{sato}, Proposition 11.10). Hence we can calculate the above integral with respect to $\nu^m$ as an integral with respect to $\nu^m_1$ and write the expression we are investigating as
\begin{equation*}
 2\mathbb{E} \int_0^t \bar{f}_{\varepsilon}(|Z_{s-}|) \left( \int_{\mathbb{R}} |y|^2  \mathbf{1}_{\{ |y| < \varepsilon/2 \}} \mathbf{1}_{\{ y(|Z_{s-}| + y)< 0 \}} \nu^m_1(dy) \right) \mathbf{1}_{\{ |Z_{s-}| > \delta \}} ds \,.
\end{equation*}
The condition $y(|Z_{s-}| + y)< 0$ holds if and only if $y < 0$ and $y \geq -|Z_{s-}|$, so for those $s \in [0,t]$ for which $|Z_{s-}| \geq \delta$ holds, we can bound the above integral with respect to $\nu^m_1$ from below, i.e.,
\begin{equation}\label{integralCepsilonDelta}
 \begin{split}
  \int_{\mathbb{R}} |y|^2 \mathbf{1}_{\{ |y| < \varepsilon/2 \}} \mathbf{1}_{\{ y(|Z_{s-}| + y)< 0 \}} \nu^m_1(dy) &\geq \int_{\mathbb{R}} |y|^2 \mathbf{1}_{\{ |y| < \varepsilon/2 \}} \mathbf{1}_{\{ y < 0 \wedge y > - \delta \}} \nu^m_1(dy) \\
  &\geq \int_{ \max \{-\delta. -\varepsilon/2 \} }^0 |y|^2 \nu^m_1(dy) > 0 \,.
 \end{split}
\end{equation}
Obviously, since $\varepsilon \leq \delta$, we have $\max \{-\delta. -\varepsilon/2 \} = -\varepsilon/2$. Moreover, we can take $m$ in our construction large enough so that
\begin{equation}\label{choiceOfm1}
 \int_{ -\varepsilon/2 }^0 |y|^2 \nu^m_1(dy) \geq \int_{ -\varepsilon/4  }^0 |y|^2 \nu_1(dy) \,,
\end{equation}
where $\nu_1$ is the first marginal of $\nu$ (note that if the dimension is greater than one, the measures $\mathbf{1}_{\{ |y| \leq m \}} \nu_1(dy)$ and $\nu^m_1(dy)$ do not coincide and hence we need to change the integration limit on the right hand side above). Thus we can estimate
\begin{equation*}
\begin{split}
 2\mathbb{E} \int_0^t  \bar{f}_{\varepsilon}(|Z_{s-}|) &\left( \int_{\mathbb{R}} |y|^2  \mathbf{1}_{\{ |y| < \varepsilon/2 \}} \mathbf{1}_{\{ y(|Z_{s-}| + y)< 0 \}} \nu_1(dy) \right) \mathbf{1}_{\{ |Z_{s-}| > \delta \}} ds \\
 &\leq C_{\varepsilon} \mathbb{E} \int_0^t \bar{f}_{\varepsilon}(|Z_{s-}|) \mathbf{1}_{\{ |Z_{s-}| > \delta \}} ds \,.
 \end{split}
\end{equation*}
\end{proof}

\end{lemma}

The calculations in the above lemma still hold if we replace the time $t$ with $t \wedge \tau_n$. Hence, after writing down the formula (\ref{ItoOptimal}) for the stopped process $(Z_{t \wedge \tau_n})_{t \geq 0}$, taking the expectation and using Lemmas \ref{lemma31}-\ref{lemma33}, we obtain
\begin{equation}\label{estimationAfterSupremum}
  \begin{split}
 \mathbb{E}f_1(|Z_{t \wedge \tau_n}|) - \mathbb{E}f_1(|Z_0|) &\leq \mathbb{E}\int_0^{t \wedge \tau_n} f_1'(|Z_{s-}|)\frac{1}{|Z_{s-}|}\langle Z_{s-} , b(X_{s-}) - b(Y_{s-}) \rangle ds \\
 &+ C_{\varepsilon} \mathbb{E} \int_0^{t \wedge \tau_n} \bar{f}_{\varepsilon}(|Z_{s-}|) \mathbf{1}_{\{ |Z_{s-}| > \delta \}} ds \,.
 \end{split}
\end{equation}
We have managed to use the second derivative of $f_1$ to obtain a negative term that works only when $|Z_{s-}| > \delta$. Recall that it was necessary to bound $|Z_{s-}|$ from below since we needed to bound the integral in (\ref{integralCepsilonDelta}) from below. In order to obtain a negative term for $|Z_{s-}| \leq \delta$ we will later use the discontinuous part $f_2$ of our distance function $f$. Now we focus on finding a continuous function $f_1$ that will give us proper estimates for $|Z_{s-}| > \delta$. The argument we use here is based on arguments used by Eberle for diffusions in his papers \cite{eberleComptes} and \cite{eberle}.

\begin{lemma}\label{lemmaFunctionalInequality}
There exist a concave, strictly increasing $\mathcal{C}^2$ function $f_1$ and a constant $c_1 > 0$ defined by (\ref{defF}) and (\ref{defC}) respectively, such that
\begin{equation}\label{functionalInequality}
 -f_1'(r)\kappa(r)r + C_{\varepsilon} \bar{f}_{\varepsilon}(r) \leq -c_1 f_1(r)
\end{equation}
holds for all $r > \delta$, where $\kappa$ is the function defined by (\ref{defKappa}).
\begin{proof}
Our assertion (\ref{functionalInequality}) is equivalent to
\begin{equation*}
 C_{\varepsilon} \bar{f}_{\varepsilon}(r) \leq -c_1 f_1(r) + f_1'(r)\kappa(r)r  \text{ for all } r > \delta
\end{equation*}
or, explicitly,
\begin{equation}\label{diffIneq1}
 \sup_{x \in (r-\varepsilon, r)} f_1''(x) \leq -\frac{c_1}{C_{\varepsilon}}f_1(r) + f_1'(r)\frac{r\kappa(r)}{C_{\varepsilon}} \text{ for all } r > \delta \,.
\end{equation}
Observe that for this to make sense, we should have $\delta \geq \varepsilon$. Define
\begin{equation*}
 h(r) := r\kappa(r) \,.
\end{equation*}
If we use the fact that $-h^{-} \leq h$, where $h^{-}$ is the negative part of $h$, then we see that in order to show (\ref{diffIneq1}), it is sufficient to show
\begin{equation*}
 \sup_{x \in (r-\varepsilon, r)} f_1''(x) \leq -\frac{c_1}{C_{\varepsilon}}f_1(r) - f_1'(r)\frac{h^{-}(r)}{C_{\varepsilon}} \text{ for all } r > \delta  \,,
\end{equation*}
which is equivalent to
\begin{equation*}
 f_1''(r-a) \leq -\frac{c_1}{C_{\varepsilon}}f_1(r) - f_1'(r)\frac{h^{-}(r)}{C_{\varepsilon}} \text{ for all } a \in (0,\varepsilon) \text{ and } r > \delta \,.
\end{equation*}
We will look for $f_1$ such that
\begin{equation*}
 f_1'(r) = \phi(r)g(r)
\end{equation*}
for some appropriately chosen functions $\phi$ and $g$. Then of course
\begin{equation*}
f_1''(r-a) = \phi'(r-a)g(r-a) + \phi(r-a)g'(r-a) \,.
\end{equation*}
We will choose $\phi$ and $g$ in such a way that
\begin{equation}\label{choiceOfG}
 \phi(r-a)g'(r-a) \leq -\frac{c_1}{C_{\varepsilon}}f_1(r)
\end{equation}
and
\begin{equation}\label{choiceOfPhi}
 \phi'(r-a)g(r-a) \leq -f_1'(r) \frac{h^{-}(r)}{C_{\varepsilon}} \,.
\end{equation}
Since we assume that $f_1'' \leq 0$, which means $f_1'$ is decreasing, we have $f_1'(r) \leq f_1'(r-a)$ and (\ref{choiceOfPhi}) is implied by
\begin{equation}\label{hereWeUseDecreasingDerivative}
 \phi'(r-a)g(r-a) \leq -f_1'(r-a) \frac{h^{-}(r)}{C_{\varepsilon}} \,.
\end{equation}
Note that our ability to replace (\ref{choiceOfPhi}) with the above condition is a consequence of our choice to consider only the jumps that decrease the distance between $X_t$ and $Y_t$ (see Remark \ref{remarkInterval}), which is equivalent to considering the supremum of $f_1''$ over a non-symmetric interval. In order to obtain (\ref{hereWeUseDecreasingDerivative}), we need $\phi$ such that
\begin{equation*}
 \phi'(r-a) \leq -\frac{h^{-}(r)}{C_{\varepsilon}}\phi(r-a) \text{ for all } a \in (0,\varepsilon) \text{ and } r > \delta \,,
\end{equation*}
which is implied by
\begin{equation}\label{choiceOfPhi2}
 \phi'(r) \leq -\frac{h^{-}(r+a)}{C_{\varepsilon}}\phi(r) \text{ for all } a \in (0,\varepsilon) \text{ and } r > 0 \,.
\end{equation}
Define
\begin{equation*}
 \bar{h}(r):=\sup_{t \in (r,r+\varepsilon)}h^{-}(t) = \sup_{t \in (r,r+\varepsilon)} t \kappa^{-}(t) \,.
\end{equation*}
Then of course
\begin{equation*}
 -\bar{h}(r) \leq -h^{-}(r+a) \text{ for all } a \in (0,\varepsilon)
\end{equation*}
and thus the condition
\begin{equation*}
 \phi'(r) \leq -\frac{\bar{h}(r)}{C_{\varepsilon}}\phi(r) \text{ for all } r > 0
\end{equation*}
implies (\ref{choiceOfPhi2}). In view of the above considerations, we can choose $\phi$ by setting
\begin{equation}\label{defPhi}
 \phi(r) := \exp{\left( -\int_0^r \frac{\bar{h}(t)}{C_{\varepsilon}} dt \right)}
\end{equation}
and this ensures that (\ref{choiceOfPhi}) holds.

If we assume $f_1(0) = 0$, then
\begin{equation}\label{defF}
 f_1(r) = \int_0^r \phi(s)g(s) ds \,.
\end{equation}
We will choose $g$ such that $1/2 \leq g \leq 1$, which will give us both a lower and an upper bound on $f_1'$. We would also like $g$ to be constant for large arguments in order to make $f_1'(r)$ constant for sufficiently large $r$. This is necessary to get an upper bound for the $W_1$ distance (see the proof of Corollary \ref{corollary2}). Hence, we will now proceed to find a formula for $g$ for which (\ref{choiceOfG}) holds and then we will extend $g$ as a constant function equal to $1/2$ beginning from some point $R_1$. Next we will show that if $R_1$ is chosen to be sufficiently large, then (\ref{diffIneq1}) holds for $r \geq R_1$ and $g = 1/2$. Note that if we set
\begin{equation*}
 \Phi(r) := \int_0^r \phi(s) ds \,,
\end{equation*}
then we have $f_1(r) \leq \Phi(r)$ and in order to get (\ref{choiceOfG}) it is sufficient to choose $g$ in such a~way that
\begin{equation}\label{choiceOfG2}
 \phi(r-a)g'(r-a) \leq -\frac{c_1}{C_{\varepsilon}}\Phi(r) \text{ for all } a \in (0, \varepsilon) \text{ and } r > \delta \,,
\end{equation}
which is implied by
\begin{equation*}
 \phi(r)g'(r) \leq - \frac{c_1}{C_{\varepsilon}} \Phi(r+a) \text{ for all } a \in (0, \varepsilon) \text{ and } r > 0 \,.
\end{equation*}
Since $\Phi$ is increasing, the condition
\begin{equation*}
 \phi(r)g'(r) \leq -\frac{c_1}{C_{\varepsilon}} \Phi(r+\varepsilon) \text{ for all } a \in (0, \varepsilon) \text{ and } r > 0
\end{equation*}
implies (\ref{choiceOfG2}). This means that we can choose $g$ by setting
\begin{equation*}
 g(r):= 1 - \frac{c_1}{C_{\varepsilon}} \int_0^r \frac{\Phi(t+\varepsilon)}{\phi(t)} dt \,.
\end{equation*}
Then obviously we have $g \leq 1$ and if we want to have $g \geq 1/2$, we must choose the constant $c_1$ in such a way that
\begin{equation*}
 1 - \frac{c_1}{C_{\varepsilon}} \int_0^r \frac{\Phi(t+\varepsilon)}{\phi(t)} dt \geq \frac{1}{2}
\end{equation*}
or equivalently
\begin{equation}\label{cInequality}
 c_1 \leq \frac{C_{\varepsilon}}{2} \left(\int_0^r \frac{\Phi(t+\varepsilon)}{\phi(t)} dt \right)^{-1} \,.
\end{equation}

Now define
\begin{equation}\label{defR0}
 R_0 := \inf \left\{R \geq 0 : \forall r \geq R : \kappa(r) \geq 0 \right\} \,.
\end{equation}
Note that $R_0$ is finite since $\lim_{r \to \infty} \kappa(r) > 0$. For all $r \geq R_0$ we have 
\begin{equation*}
  h^{-}(r) = 0 \text{  and } \phi(r) = \phi(R_0) \,.
\end{equation*}
Now we would like to define $R_1 \geq R_0 + \varepsilon$ in such a way that
\begin{equation*}
 g(r) = \begin{cases} 
      1 - \frac{c_1}{C_{\varepsilon}} \int_0^r \frac{\Phi(t+\varepsilon)}{\phi(t)} dt & r\leq R_1 \\
      \frac{1}{2} & r\geq R_1
   \end{cases}
\end{equation*}
and (\ref{diffIneq1}) holds for $r \geq R_1$. By setting
\begin{equation}\label{defC}
 c_1 := \frac{C_{\varepsilon}}{2}\left(\int_0^{R_1} \frac{\Phi(t+\varepsilon)}{\phi(t)} dt\right)^{-1}
\end{equation}
we ensure that $g$ defined above is continuous and that (\ref{cInequality}) and, in consequence, (\ref{choiceOfG}) holds for $r \leq R_1$.

We will now explain how to find $R_1$. Since $f_1'(r) = \frac{1}{2}\phi(R_0)$ for $r \geq R_1$, we have
\begin{equation*}
 \sup_{x \in (r - \varepsilon,r)} f_1''(x) = 0 \text{ for all } r \geq R_1
\end{equation*}
and therefore (\ref{diffIneq1}) for $r \geq R_1$ holds if and only if
\begin{equation*}
 - f_1'(r)\frac{r\kappa(r)}{C_{\varepsilon}} \leq -\frac{c_1}{C_{\varepsilon}}f_1(r) \text{ for all } r \geq R_1 \,,
\end{equation*}
which is equivalent to
\begin{equation*}
 - r\kappa(r) \frac{\phi(R_0)}{2} \leq -c_1 f_1(r)  \text{ for all } r \geq R_1 \,.
\end{equation*}
Using once again the fact that $f_1 \leq \Phi$, we see that it is sufficient to have
\begin{equation*}
 - r\kappa(r) \frac{\phi(R_0)}{2} \leq -c_1 \Phi(r)  \text{ for all } r \geq R_1 \,.
\end{equation*}
By the definition of $c_1$, the right hand side of the above inequality is equal to
\begin{equation*}
 -C_{\varepsilon}\Phi(r)  \left( 2\int_0^{R_1} \frac{\Phi(t+\varepsilon)}{\phi(t)}dt \right)^{-1} \,.
\end{equation*}
In order to make our computations easier, we will use the inequality
\begin{equation*}
 \int_{R_0}^{R_1} \frac{\Phi(t+\varepsilon)}{\phi(t)}dt \leq \int_0^{R_1} \frac{\Phi(t+\varepsilon)}{\phi(t)}dt
\end{equation*}
and we will look for $R_1$ such that
\begin{equation}\label{lookingForR}
 - r\kappa(r) \frac{\phi(R_0)}{2} \leq -C_{\varepsilon}\Phi(r) \left( 2\int_{R_0}^{R_1} \frac{\Phi(t+\varepsilon)}{\phi(t)}dt \right)^{-1}  \text{ for all } r \geq R_1 \,.
\end{equation}
We can compute
\begin{equation*}
 \begin{split}
  \int_{R_0}^{R_1} \frac{\Phi(t+\varepsilon)}{\phi(t)}dt &= \int_{R_0}^{R_1} \frac{\Phi(R_0)+\phi(R_0)(t+\varepsilon-R_0)}{\phi(R_0)}dt \\
  &=(R_1 - R_0)\frac{\Phi(R_0)}{\phi(R_0)} + \frac{1}{2}(R_1+\varepsilon - R_0)^2 - \frac{1}{2}\varepsilon^2 \\
  &\geq (R_1 - R_0)\frac{\Phi(R_0)}{\phi(R_0)} + \frac{1}{2}(R_1 - R_0)^2 \\
  &\geq \frac{1}{2}(R_1 - R_0)\frac{\Phi(R_0)}{\phi(R_0)} + \frac{1}{2}(R_1 - R_0)^2 \\
  &= \frac{(R_1-R_0)\Phi(R_1)}{2\phi(R_0)} \,.
 \end{split}
\end{equation*}
Therefore if we find $R_1$ such that
\begin{equation}\label{R1ineq1}
 - r\kappa(r) \frac{\phi(R_0)}{2} \leq \frac{-C_{\varepsilon}\Phi(r)\phi(R_0)}{(R_1-R_0)\Phi(R_1)}  \text{ for all } r \geq R_1 \,,
\end{equation}
it will imply (\ref{lookingForR}). Observe now that we have 
\begin{equation}\label{R1ineq2}
 \frac{\Phi(r)}{\Phi(R_1)} \leq \frac{r}{R_1} \text{ for all } r \geq R_1 \,.
\end{equation}
This follows from the fact that $\phi$ is decreasing, which implies that $\Phi(R_1) \geq \phi(R_0)R_1$ and thus
\begin{equation*}
 \frac{\phi(R_0)}{\Phi(R_1)}(r - R_1) \leq \frac{1}{R_1}(r - R_1)
\end{equation*}
and
\begin{equation*}
 \frac{\phi(R_0)(r - R_1) + \Phi(R_1)}{\Phi(R_1)} \leq \frac{r}{R_1}
\end{equation*}
hold for $r \geq R_1$. If we divide both sides of (\ref{R1ineq1}) by $\phi(R_0)$ and use (\ref{R1ineq2}), we see that we need to have
\begin{equation*}
  \frac{-r\kappa(r)}{2} \leq \frac{-C_{\varepsilon}r}{(R_1-R_0)R_1}  \text{ for all } r \geq R_1
\end{equation*}
or, equivalently,
\begin{equation*}
  \frac{2C_{\varepsilon}}{(R_1 - R_0)R_1} \leq \kappa(r)  \text{ for all } r \geq R_1 \,.
\end{equation*}
This shows that we can define $R_1$ by
\begin{equation}\label{defR1}
 R_1 := \inf \left\{ R \geq R_0 + \varepsilon : \forall r \geq R : \kappa(r) \geq \frac{2C_{\varepsilon}}{(R-R_0)R} \right\} \,,
\end{equation}
which is finite since we assume that $\lim_{r \to \infty} \kappa(r) > 0$. 
\end{proof}
\end{lemma}

Our choice of $f_1$ and $c_1$ made above (see (\ref{defF}) and (\ref{defC}), respectively) allows us to estimate
\begin{equation}\label{OutsideBall1}
 \begin{split}
 &\mathbb{E}\int_0^{t \wedge \tau_n} f_1'(|Z_{s-}|)\frac{1}{|Z_{s-}|}\langle Z_{s-} , b(X_{s-}) - b(Y_{s-}) \rangle \mathbf{1}_{\{ |Z_{s-}| > \delta \}}ds \\
 &+C_{\varepsilon}\mathbb{E} \int_0^{t \wedge \tau_n} \bar{f}_{\varepsilon}(|Z_{s-}|)\mathbf{1}_{\{ |Z_{s-}| > \delta \}} ds \leq \mathbb{E} \int_0^{t \wedge \tau_n} -c_1 f_1(|Z_s|) \mathbf{1}_{\{ |Z_{s-}| > \delta \}} ds \,.
 \end{split}
\end{equation}
If we are to obtain (\ref{preGronwall}), then on the right hand side of (\ref{OutsideBall1}) we would like to have the function $f$ instead of $f_1$, but we can achieve this by assuming
\begin{equation}\label{BoundFora}
 a \leq K \inf_{x > \delta} f_1(x)
\end{equation}
or, more explicitly, $a \leq K f_1(\delta)$ (since $f_1$ is increasing), for some constant $K \geq 1$ to be chosen later. Then we have
\begin{equation*}
\begin{split}
 -c_1 f_1(|Z_s|)\mathbf{1}_{\{ |Z_{s-}| > \delta \}} &= -c_1 \left[\frac{1}{2}f_1(|Z_s|) + \frac{1}{2}f_1(|Z_s|)\right] \mathbf{1}_{\{ |Z_{s-}| > \delta \}} \\
 &\leq -\frac{c_1}{2}f_1(|Z_s|)\mathbf{1}_{\{ |Z_{s-}| > \delta \}} -\frac{c_1 a}{2K}\mathbf{1}_{\{ |Z_{s-}| > \delta \}}\\
 &\leq -\frac{c_1}{2K} (f_1 + a)(|Z_s|)\mathbf{1}_{\{ |Z_{s-}| > \delta \}} = -\frac{c_1}{2K} f(|Z_s|)\mathbf{1}_{\{ |Z_{s-}| > \delta \}} \,.
\end{split}
\end{equation*}
and hence
\begin{equation}\label{OutsideBall2}
 \mathbb{E}\int_0^{t \wedge \tau_n} -c_1 f_1(|Z_s|)\mathbf{1}_{\{ |Z_{s-}| > \delta \}} ds \leq \mathbb{E}\int_0^{t \wedge \tau_n} -\frac{c_1}{2K} f(|Z_s|)\mathbf{1}_{\{ |Z_{s-}| > \delta \}}ds
\end{equation}
Now if we write (\ref{estimationAfterSupremum}) as
\begin{equation}\label{EstimationTwoParts}
  \begin{split}
 \mathbb{E}f_1(|Z_{t \wedge \tau_n}|) &- \mathbb{E}f_1(|Z_0|) \\
 &\leq \mathbb{E}\int_0^{t \wedge \tau_n} f_1'(|Z_{s-}|)\frac{1}{|Z_{s-}|}\langle Z_{s-} , b(X_{s-}) - b(Y_{s-}) \rangle (\mathbf{1}_{\{ |Z_{s-}| > \delta \}} + \mathbf{1}_{\{ 0 < |Z_{s-}| \leq \delta \}})ds \\
 &+ C_{\varepsilon} \mathbb{E} \int_0^{t \wedge \tau_n} \bar{f}_{\varepsilon}(|Z_{s-}|) \mathbf{1}_{\{ |Z_{s-}| > \delta \}} ds \,,
 \end{split}
\end{equation}
we see that by (\ref{OutsideBall1}) and (\ref{OutsideBall2}) we already have a good bound for the terms involving $\mathbf{1}_{\{ |Z_{s-}| > \delta \}}$. Now we need to obtain estimates for the case when $|Z_{s-}| \leq \delta$. To this end, we should come back to the equation (\ref{ItoPreliminary}) and focus on the expression 
\begin{equation*}
 a\mathbb{E}\mathbf{1}_{(0,\infty)}(|Z_t|) - a\mathbb{E}\mathbf{1}_{(0,\infty)}(|Z_0|) \,.
\end{equation*}
We have the following lemma.

\begin{lemma}\label{lemmaDiscontinuous}
 For any $t \geq 0$ we have
 \begin{equation*}
 \mathbb{E}\mathbf{1}_{(0,\infty)}(|Z_t|) - \mathbb{E}\mathbf{1}_{(0,\infty)}(|Z_0|) \leq - \mathbb{E} \int_0^t \widetilde{C}_{\delta}(m) \mathbf{1}_{\{ 0< |Z_{s-}| \leq \delta \}} ds \,,
 \end{equation*}
 where 
 \begin{equation}\label{Cdelta}
  \widetilde{C}_{\delta}(m) := \inf_{x \in \mathbb{R}^d : 0 < |x| \leq \delta} \int_{\{ |v| \leq m \} \cap \{ |v + x| \leq m \}} q(v) \wedge q(v + x) dv > 0 \,.
 \end{equation}
 Note that $\widetilde{C}_{\delta}(m)$ is positive by Assumption \ref{Assumption3} about the sufficient overlap of $q$ and translated $q$ (see the condition (\ref{overlapCondition})).
 \begin{proof}
  Observe that almost surely we have
\begin{equation*}
 \mathbf{1}_{(0,\infty)}(|Z_t|) = 1 - \int_0^t \int_{\{ |v| \leq m \} \times [0,1]} \mathbf{1}_{\{ u < \rho(v, Z_{s-}) \}} \mathbf{1}_{\{ |Z_{s-}| \neq 0 \}} N(ds,dv,du) \,.
\end{equation*}
The integral with respect to the Poisson random measure $N$ appearing above counts exactly the one jump that brings the processes $X_t$ and $Y_t$ to the same point. Note that if we skipped the condition $\{ |Z_{s-}| \neq 0 \}$, it would also count all the jumps that happen after the coupling time and it would be possibly infinite. Since we obviously have
\begin{equation*}
\begin{split}
 \int_0^t \int_{\{ |v| \leq m \} \times [0,1]} &\mathbf{1}_{\{ u < \rho(v, Z_{s-}) \}} \mathbf{1}_{\{ 0 < |Z_{s-}| \leq \delta \}} N(ds,dv,du) \\
 &\leq \int_0^t \int_{\{ |v| \leq m \} \times [0,1]} \mathbf{1}_{\{ u < \rho(v, Z_{s-}) \}} \mathbf{1}_{\{ |Z_{s-}| \neq 0 \}} N(ds,dv,du) \,,
 \end{split}
\end{equation*}
we can estimate
\begin{equation*}
  \mathbf{1}_{(0,\infty)}(|Z_t|) \leq 1 - \int_0^t \int_{\{ |v| \leq m \} \times [0,1]} \mathbf{1}_{\{ u < \rho(v, Z_{s-}) \}} \mathbf{1}_{\{ 0 < |Z_{s-}| \leq \delta \}} N(ds,dv,du) \,,
\end{equation*}
and therefore we get
\begin{equation*}
 a\mathbb{E}\mathbf{1}_{(0,\infty)}(|Z_t|) -  a\mathbb{E}\mathbf{1}_{(0,\infty)}(|Z_0|) \leq -a \mathbb{E} \int_0^t \int_{\{ |v| \leq m \} \times [0,1]} \mathbf{1}_{\{ u < \rho(v, Z_{s-}) \}} \mathbf{1}_{\{ 0 < |Z_{s-}| \leq \delta \}} \nu(dv) du ds \,,
\end{equation*}
where we used the assumption that $\mathbb{E} \mathbf{1}_{(0,\infty)}(|Z_0|) = \mathbb{P}(|Z_0| \neq 0) = 1$ (see the remarks at the beginning of this section). We also have
\begin{equation*}
\begin{split}
 \mathbb{E} \int_0^t \int_{\{ |v| \leq m \} \times [0,1]} &\mathbf{1}_{\{ u < \rho(v, Z_{s-}) \}} \mathbf{1}_{\{ 0 < |Z_{s-}| \leq \delta \}} \nu(dv) du ds \\
 &= \mathbb{E} \int_0^t \int_{\{ |v| \leq m \}} \rho(v, Z_{s-}) \mathbf{1}_{\{ 0< |Z_{s-}| \leq \delta \}} \nu(dv) ds\\
 &= \mathbb{E} \int_0^t \int_{\{ |v| \leq m \} \cap \{ |v + Z_{s-}| \leq m \}} (q(v + Z_{s-}) \wedge q(v)) \mathbf{1}_{\{ 0< |Z_{s-}| \leq \delta \}} dv ds \\
 &\geq \mathbb{E} \int_0^t \widetilde{C}_{\delta}(m) \mathbf{1}_{\{ 0< |Z_{s-}| \leq \delta \}} ds
\end{split}
\end{equation*}
and the assertion follows.
 \end{proof}

\end{lemma}

Note that we can always choose $m$ large enough so that
\begin{equation}\label{choiceOfm2}
 \inf_{x \in \mathbb{R}^d : 0 < |x| \leq \delta} \int_{\{ |v| \leq m \} \cap \{ |v + x| \leq m \}} q(v) \wedge q(v + x) dv \geq \frac{1}{2} \inf_{x \in \mathbb{R}^d : 0 < |x| \leq \delta} \int_{\mathbb{R}^d} q(v) \wedge q(v + x) dv =: \frac{1}{2} \widetilde{C}_{\delta}
\end{equation}
and hence we have
\begin{equation*}
 \mathbb{E}\mathbf{1}_{(0,\infty)}(|Z_t|) - \mathbb{E}\mathbf{1}_{(0,\infty)}(|Z_0|) \leq - \mathbb{E} \int_0^t \frac{1}{2} \widetilde{C}_{\delta} \mathbf{1}_{\{ 0< |Z_{s-}| \leq \delta \}} ds \,.
\end{equation*}
Combining the estimate above with (\ref{ItoPreliminary}) and (\ref{EstimationTwoParts}), we obtain
\begin{equation*}
 \begin{split}
 \mathbb{E}f(|Z_{t \wedge \tau_n}|) &- \mathbb{E}f(|Z_0|) \\
 &\leq \mathbb{E}\int_0^{t \wedge \tau_n} f_1'(|Z_{s-}|)\frac{1}{|Z_{s-}|}\langle Z_{s-} , b(X_{s-}) - b(Y_{s-}) \rangle (\mathbf{1}_{\{ |Z_{s-}| > \delta \}} + \mathbf{1}_{\{ 0 < |Z_{s-}| \leq \delta \}})ds \\
 &+ C_{\varepsilon} \mathbb{E} \int_0^{t \wedge \tau_n} \bar{f}_{\varepsilon}(|Z_{s-}|) \mathbf{1}_{\{ |Z_{s-}| > \delta \}} ds - a\mathbb{E} \int_0^{t \wedge \tau_n} \frac{1}{2} \widetilde{C}_{\delta} \mathbf{1}_{\{ 0< |Z_{s-}| \leq \delta \}} ds \,.
 \end{split}
\end{equation*}
In order to deal with the expressions involving $\{ |Z_{s-}| \leq \delta \}$, we will use the fact that $b$ satisfies the one-sided Lipschitz condition (\ref{onesidedLipschitz}) with some constant $C_L > 0$ and that $f_1'(r) \leq f_1'(0) = 1$ for all $r \geq 0$ to get
\begin{equation}\label{ExpressionToBound}
 \begin{split}
 \mathbb{E}\int_0^{t \wedge \tau_n} &f_1'(|Z_{s-}|) \frac{1}{|Z_{s-}|}\langle Z_{s-} , b(X_{s-}) - b(Y_{s-}) \rangle \mathbf{1}_{\{ 0 < |Z_{s-}| \leq \delta \}}ds - a\mathbb{E} \int_0^{t \wedge \tau_n} \frac{1}{2} \widetilde{C}_{\delta} \mathbf{1}_{\{ 0< |Z_{s-}| \leq \delta \}} ds \\
 &\leq (C_L \delta - \frac{1}{2}a\widetilde{C}_{\delta}) \mathbb{E}\int_0^{t \wedge \tau_n} \mathbf{1}_{\{ 0 < |Z_{s-}| \leq \delta \}}ds \,.
 \end{split}
\end{equation}
We would like to bound this expression by
\begin{equation*}
 \mathbb{E}\int_0^{t \wedge \tau_n} -Cf(|Z_{s-}|)\mathbf{1}_{\{ 0 < |Z_{s-}| \leq \delta \}} ds
\end{equation*}
for some positive constant $C$, but since the function $f$ is bounded on the interval $[0,\delta]$ by $f_1(\delta) + a$, we have
\begin{equation*}
 -Cf_1(\delta) - Ca \leq -Cf(|Z_{s-}|) \text{ if } 0 < |Z_{s-}| \leq \delta
\end{equation*}
and thus it is sufficient if we have
\begin{equation*}
 C_L\delta + Cf_1(\delta) \leq ( \widetilde{C}_{\delta}/2 - C)a \,.
\end{equation*}
Of course the right hand side has to be positive, so we can choose e.g. $C := \widetilde{C}_{\delta} / 4$.
Then we must have
\begin{equation}\label{LowerBoundFora}
 \frac{C_L\delta + \widetilde{C}_{\delta}f_1(\delta)/4}{\widetilde{C}_{\delta} / 4} \leq a \,,
\end{equation}
but on the other hand, by (\ref{BoundFora}), we must also have $a \leq Kf_1(\delta)$.
Hence we can define
\begin{equation}\label{defK}
 K :=  \frac{C_L\delta + \widetilde{C}_{\delta}f_1(\delta)/4}{\widetilde{C}_{\delta}f_1(\delta)/4}
\end{equation}
and
\begin{equation}\label{defA}
 a:=Kf_1(\delta) \,.
\end{equation}
Then obviously both (\ref{LowerBoundFora}) and (\ref{BoundFora}) hold and we get the required estimate for the right hand side of (\ref{ExpressionToBound}). Using all our estimates together, we get
\begin{equation*}
\begin{split}
 \mathbb{E}f(|Z_{t \wedge \tau_n}|) - \mathbb{E}f(|Z_0|) &\leq \mathbb{E} \int_0^{t \wedge \tau_n} -\frac{c_1}{2K}f(|Z_s|)\mathbf{1}_{\{ |Z_{s-}| > \delta \}} ds \\
 &+ \mathbb{E} \int_0^{t \wedge \tau_n} -\frac{1}{4}\widetilde{C}_{\delta} f(|Z_s|)\mathbf{1}_{\{ 0 < |Z_{s-}| \leq \delta \}} ds \,.
 \end{split}
\end{equation*}
Denote 
\begin{equation}\label{contractivityConstant}
c:=\min \left\{\frac{c_1}{2K}, \frac{1}{4}\widetilde{C}_{\delta}\right\} \,.
\end{equation}
Then of course
\begin{equation}\label{almostFinalInequality}
\begin{split}
 \mathbb{E}f(|Z_{t \wedge \tau_n}|) - \mathbb{E}f(|Z_0|) &\leq \mathbb{E} \int_0^{t \wedge \tau_n} -cf(|Z_s|)\mathbf{1}_{\{ |Z_{s-}| > \delta \}} ds \\
 &+ \mathbb{E} \int_0^{t \wedge \tau_n} -cf(|Z_s|)\mathbf{1}_{\{ 0 < |Z_{s-}| \leq \delta \}} ds \\
 &= \mathbb{E} \int_0^{t \wedge \tau_n} -cf(|Z_s|) ds \,.
\end{split}
\end{equation}
Note that we can perform the same calculations not only on the interval $[0,t \wedge \tau_n]$, but also on any interval $[s \wedge \tau_n, t \wedge \tau_n]$ for arbitrary $0 \leq s < t$. Indeed, by our assumption (see the beginning of this section) we have $\mathbb{P}(|Z_0| \neq 0) = 1$ and hence for any $0 \leq s < T$ we have $\mathbb{P}(|Z_{s \wedge \tau_n}| \neq 0) = 1$. Thus Lemma \ref{lemmaDiscontinuous} still holds on $[s \wedge \tau_n, t \wedge \tau_n]$. It is easy to see that the other calculations are valid too and we obtain
\begin{equation}\label{preDiffGronwall}
 \mathbb{E}f(|Z_{t \wedge \tau_n}|) - \mathbb{E}f(|Z_{s \wedge \tau_n}|) \leq \mathbb{E} \int_s^t -c f(|Z_{r \wedge \tau_n}|) dr \,.
\end{equation}
Since this holds for any $0 \leq s < t$, by the differential version of the Gronwall inequality we obtain
\begin{equation*}
 \mathbb{E}f(|Z_{t \wedge \tau_n}|) \leq \mathbb{E}f(|Z_0|)e^{-ct} \,.
\end{equation*}
Note that we cannot use the integral version of the Gronwall inequality for (\ref{almostFinalInequality}) since the right hand side is negative and that is why we need (\ref{preDiffGronwall}) to hold for any $s < t$. By the Fatou lemma and the fact that $Z_t = 0$ for $t \geq T$ (see the remarks after (\ref{ItoOptimal})) we get
\begin{equation*}
 \mathbb{E}f(|Z_t|) \leq e^{-ct} \mathbb{E}f(|Z_0|) \text{ for all } t \geq 0 \,,
\end{equation*}
which finishes the proof of (\ref{mainInequality}).

\begin{proof0}
By everything we proved in Sections \ref{sectionExistence}, \ref{sectionCoupling} and the entire Section \ref{Section3}, we obtain a coupling $(X_t, Y_t)_{t \geq 0}$ satisfying the inequality (\ref{mainInequality}). The only thing that remains to be shown is the fact that the coupling $(X_t, Y_t)_{t \geq 0}$ is successful. This follows easily from the inequality (\ref{mainInequality}) and the form of the function $f$. Indeed, recalling that $Z_t = X_t - Y_t$ and that $T$ denotes the coupling time for $(X_t, Y_t)_{t \geq 0}$, for a fixed $t > 0$ we have
\begin{equation*}
\begin{split}
 \mathbb{P}(T > t) &= \mathbb{P}(|Z_t|>0) = \mathbb{E} \mathbf{1}_{(0,\infty)}(|Z_t|) \leq \frac{1}{a} \mathbb{E} \left( f_1(|Z_t|) + a \mathbf{1}_{(0,\infty)}(|Z_t|) \right) \\
 &= \frac{1}{a} \mathbb{E} f(|Z_t|) \leq \frac{1}{a}e^{-ct} \mathbb{E}f(|Z_0|) \,.
 \end{split}
\end{equation*}
Hence we get
\begin{equation*}
 \mathbb{P}(T = \infty) = \mathbb{P} \left( \bigcap_{t > 0} \{ T > t \} \right) = \lim_{t \to \infty} \mathbb{P}(T > t) = 0 \,.
\end{equation*}
 \qed
\end{proof0}

\section{Additional proofs and examples}\label{Section4}

\begin{proof1}
We have 
\begin{equation*}
 \mathbf{1}_{(0,\infty)} = a^{-1}a\mathbf{1}_{(0,\infty)} \leq a^{-1}(f_1 + a\mathbf{1}_{(0,\infty)}) = a^{-1}f \,,
\end{equation*}
hence we get
\begin{equation*}
 \frac{1}{2}\| \mu_1 p_t - \mu_2 p_t \|_{TV} = W_{\mathbf{1}_{(0,\infty)}}(\mu_1 p_t, \mu_2 p_t) \leq a^{-1}W_f(\mu_1 p_t, \mu_2 p_t) \leq a^{-1} e^{-ct} W_f(\mu_1, \mu_2) \,.
\end{equation*}
\qed
\end{proof1}

\begin{proof2}
We have
\begin{equation*}
 f_1'(r) = \phi(r)g(r) \geq \frac{\phi(r)}{2} \geq \frac{\phi(R_0)}{2}
\end{equation*}
for all $r \geq 0$. But $f_1(0) = 0$, so we get
\begin{equation*}
 f_1(r) \geq \frac{\phi(R_0)}{2}r
\end{equation*}
for all $r \geq 0$ and in consequence
\begin{equation*}
 r \leq \frac{2f_1(r)}{\phi(R_0)} \leq \frac{2f(r)}{\phi(R_0)} \,,
\end{equation*}
which proves that
\begin{equation*}
 W_1(\mu_1 p_t , \mu_2 p_t) \leq 2\phi(R_0)^{-1} e^{-ct} W_f(\mu_1, \mu_2) \,.
\end{equation*}
\qed
\end{proof2}

\begin{proof3}
Let us first comment on the assumption we make on the semigroup $(p_t)_{t \geq 0}$ stating that if a measure $\mu$ has a finite first moment, then for all $t > 0$ the measure $\mu p_t$ also has a finite first moment. This assumption seems quite natural for proving existence of invariant measures for Markov processes by using methods based on Wasserstein distances, cf. assumption $(H1)$ in \cite{komorowski}. In our setup, it holds e.g. if we assume that the noise $(L_t)_{t \geq 0}$ has a finite first moment and the drift $b$ satisfies a linear growth condition, i.e., there exists a constant $C > 0$ such that $|b(x)|^2 \leq C(1 + |x^2|)$ for all $x \in \mathbb{R}^d$.
By Corollary \ref{corollary2}, we have
\begin{equation}\label{W1upperbound2}
 W_1(\mu p_t , \eta p_t) \leq L e^{-ct} W_f(\mu, \eta)
\end{equation}
for some constants $c$, $L > 0$ and any probability measures $\mu$ and $\eta$. Now let $\mu$ be a fixed, arbitrarily chosen probability measure and consider a sequence of measures $(\mu p_n)_{n = 0}^{\infty}$. Apply (\ref{W1upperbound2}) to $\mu$ and $\eta = \mu p_1$ with $t = n$. We get
\begin{equation*}
  W_1(\mu p_n , \mu p_{n+1}) \leq L e^{-cn} W_f(\mu, \mu p_1) \,.
\end{equation*}
Similarly, using the triangle inequality for $W_1$, we get that for any $k \geq 1$
\begin{equation*}
W_1(\mu p_n , \mu p_{n+k}) \leq L \sum_{j=0}^{k-1} e^{-c(n+j)} W_f(\mu, \mu p_1) \leq L \frac{e^{-cn}}{1 - e^{-c}} W_f(\mu, \mu p_1) \,.
\end{equation*}
It is now easy to see that $(\mu p_n)_{n = 0}^{\infty}$ is a Cauchy sequence with respect to the $W_1$ distance. Since the space of probability measures with finite first moments equipped with the $W_1$ distance is complete (see e.g. Theorem 6.18 in \cite{villani}), we infer that $(\mu p_n)_{n = 0}^{\infty}$ has a limit $\mu_0$. Note that here we use the assumption about the semigroup $(p_t)_{t \geq 0}$ preserving finite first moments. 
We also know that $W_1$ actually metrizes the weak convergence of measures and thus
\begin{equation*}
 \int \varphi \mu p_n \to \int \varphi \mu_0
\end{equation*}
as $n \to \infty$ for all continuous bounded ($\mathcal{C}_b$) functions $\varphi$. It is easy to check that since the drift in (\ref{SDE1}) is one-sided Lipschitz, the semigroup $(p_t)_{t \geq 0}$ is Feller, in particular for any $\varphi \in \mathcal{C}_b$ we have $p_1 \varphi \in \mathcal{C}_b$ and thus
\begin{equation*}
  \int \varphi(x) \mu p_{n+1}(dx) = \int p_1 \varphi (x) \mu p_n(dx) \to \int p_1 \varphi (x) \mu_0 (dx) = \int \varphi (x) \mu_0 p_1 (dx) \,.
\end{equation*}
Hence we infer that
\begin{equation*}
 \mu_0 = \mu_0 p_1 \,.
\end{equation*}
Now if we define
\begin{equation*}
 \mu_{*} := \int_0^1 \mu_0 p_s ds \,,
\end{equation*}
we can easily show (see e.g. \cite{komorowski}, the beginning of Section 3 for details) that for any $t \geq 0$ we have
\begin{equation*}
 \mu_{*} p_t = \mu_{*} \,,
\end{equation*}
i.e., $\mu_{*}$ is actually an invariant measure for $(p_t)_{t \geq 0}$. 
Now the inequality (\ref{Wfinvariantbound}) follows easily from (\ref{Wfbound}) applied to $\mu_{*}$ and $\eta$. Indeed, we have
\begin{equation*}
 W_f(\mu_{*}, \eta p_t) = W_f(\mu_{*} p_t , \eta p_t) \leq e^{-ct} W_f(\mu_{*}, \eta) \,.
\end{equation*}
Similarly, the inequalities (\ref{TVinvariantbound}) and (\ref{W1invariantbound}) follow easily from (\ref{TVupperbound}) and (\ref{W1upperbound}), respectively.
\qed
\end{proof3}

We would like now to investigate optimality of the contraction constant we obtained in Corollary \ref{mainTheorem}. First, let us recall a well-known result. Let $(X_t)_{t \geq 0}$ be the solution to (\ref{SDE1}) and $(p_t)_{t \geq 0}$ its associated semigroup. If there exists a constant $M > 0$ such that for all $x$, $y \in \mathbb{R}^d$ we have
 \begin{equation}\label{convexityAssumption}
  \langle b(x) - b(y) , x - y \rangle \leq -M|x-y|^2 \,,
 \end{equation}
then for all $t > 0$ and any probability measures $\mu_1$, $\mu_2$ we have
\begin{equation*}
 W_1(\mu_1 p_t , \mu_2 p_t) \leq e^{-Mt} W_1(\mu_1 , \mu_2) \,.
\end{equation*}

\begin{example}
A typical example illustrating the above result is the case when the drift $b$ is given as the gradient of a convex potential, i.e., $b = - \nabla U$ with e.g. $U(x) = M|x^2|/2$ for some constant $M >0$. Then we obviously have
\begin{equation*}
 \langle b(x) - b(y) , x - y \rangle = -M|x-y|^2
\end{equation*}
and, by the above result, exponential convergence with the rate $e^{-Mt}$ holds for the equation (\ref{SDE1}) in the standard $L^1$-Wasserstein distance.
\end{example}

\begin{example}\label{exampleConstantBounds}
We will now try to examine the case in which we drop the convexity assumption.
Assume 
\begin{equation}\label{nonConvexCase}
 \kappa(r) \geq 0 \text{  for all } r \geq 0 \text{  and  } \kappa(r) \geq M \text{  for all } r \geq R
\end{equation}
for some constants $M > 0$ and $R > 0$. This means that we have
\begin{equation*}
\langle b(x) - b(y) , x - y \rangle \leq 0 
\end{equation*}
everywhere, but the dissipativity condition (\ref{convexityAssumption}) holds only outside some fixed ball of radius $R$. 
Then, using the notation from Section \ref{Section3}, we can easily check that the function $\phi$ is constant and equal to $1$. We have
\begin{equation*}
 f_1(r) = \int_0^r g(s) ds \text{  and }  g(r) = 1 -  \frac{1}{R_1^2 + 2\varepsilon R_1}\left(\frac{1}{2}r^2 + \varepsilon r\right)
\end{equation*}
and therefore
\begin{equation*}
 f_1(r) = r - \frac{1}{R_1^2 + 2\varepsilon R_1}\left(\frac{1}{6}r^3 + \frac{1}{2}\varepsilon r^2\right) \,.
\end{equation*}
We also have $R_0 = 0$ and it can be shown that
\begin{equation*}
  R_1 \leq \max (R, W) \,,
\end{equation*}
where $W$ is the positive solution to the equation $M = 2C_{\varepsilon} / W^2 $, i.e., $W = \sqrt{2C_{\varepsilon} / M}$.
Indeed, if $R > W$, then $2C_{\varepsilon} / R^2 \leq 2C_{\varepsilon} / W^2 = M$ and thus, by (\ref{nonConvexCase}), for all $r \geq R$ we have $\kappa(r) \geq 2C_{\varepsilon} / R^2$, which implies that $R$ belongs to the set of which $R_1$ is the infimum (see (\ref{defR1})) and hence $R_1 \leq R$. On the other hand, if $R \leq W$, then for all $r \geq W$ we have $\kappa(r) \geq M = 2C_{\varepsilon} / W^2$ and thus $R_1 \leq W$. Observe that
\begin{equation*}
 c_1 = \frac{C_{\varepsilon}}{ R_1^2 + 2\varepsilon R_1} \geq \frac{C_{\varepsilon}}{ \max (R, W)^2 + 2 \varepsilon \max (R, W)} \,.
\end{equation*}
Moreover, $K = 1$ when $C_L = 0$ (see (\ref{defK})). Thus we have
\begin{equation*}
 \frac{c_1}{2K} \geq  \frac{C_{\varepsilon}}{ 2 \max (R, W)^2 + 4 \varepsilon \max (R, W)} \,,
\end{equation*}
which means that the lower bound for $c_1 / 2K$ is of order $\min (R^{-2},M)$. This means that the convergence rates in the $W_1$ distance are not substantially affected by dropping the global dissipativity assumption, as long as the ball in which the dissipativity does not hold is not too large. This behaviour is similar to the diffusion case (see Remark 5 in \cite{eberle}).

As an example, consider a one-dimensional L\'{e}vy process with the jump density given by $q(v) = (1/|v|^{1+\alpha})$ for $\alpha \in (0,2)$. Then we can easily show that 
\begin{equation*}
 C_{\varepsilon} = \frac{2}{2-\alpha}\left(\frac{\varepsilon}{4}\right)^{2 - \alpha} \text{ and  } \widetilde{C}_{\delta} = \frac{2}{\alpha}\left(\frac{2}{\delta}\right)^{\alpha} \,.
\end{equation*}
Let us focus on the case of $\alpha \in (1,2)$. If we denote
\begin{equation*}
 c_1(\varepsilon) :=  \frac{C_{\varepsilon}}{ 2 R^2 + 4 \varepsilon R} \,,
\end{equation*}
then as a function of $\varepsilon$ it obtains its maximum for $\varepsilon_0 := (2 - \alpha)R(2\alpha - 2)^{-1}$. Thus if $c_1(\varepsilon_0) \leq c_2(\varepsilon_0)$, where $c_2(\delta) := \widetilde{C}_{\delta}/4$ (which, as we can check numerically, is true e.g. for any $R$ if $\alpha > 11/10$), then we see that the optimal choice of parameters that maximizes the lower bound for $c = \min \{ c_1 / 2K , \widetilde{C}_{\delta}/4 \}$ is to take $\varepsilon = \delta = \varepsilon_0$, at least as long as $R \geq \sqrt{2C_{\varepsilon_0}/M}$, since only then $c_1(\varepsilon_0)$ is actually a lower bound for $c_1/2K$. But for this to be true, once $R$ and $\alpha$ are fixed, it is sufficient to consider a large enough $M$ (to give specific values, e.g. for $R=1$ and $\alpha = 3/2$ we have $\varepsilon_0 = 1/2$, $C_{\varepsilon_0} = \sqrt{2}$ and $c_1(\varepsilon_0) = \sqrt{2}/4$, hence when we consider $M \geq 2 \sqrt{2}$, it is optimal to take $\varepsilon = \delta = 1/2$ and we obtain $c \geq \sqrt{2}/8$). Note that for fixed values of $R$ and $M$, when $\alpha$ increases to $2$, the values of $C_{\varepsilon_0}$, $c_1(\varepsilon_0)$ and $c_2(\varepsilon_0)$ increase to $\infty$. However, in such a case $c_1(\varepsilon_0)$ is no longer a lower bound for $c_1/2K$, since $R < \sqrt{2C_{\varepsilon_0}/M}$. Instead we have
\begin{equation*}
 \frac{c_1}{2K} \geq \frac{C_{\varepsilon_0}}{4C_{\varepsilon_0}M^{-1} + 4\varepsilon_0 \sqrt{2 C_{\varepsilon_0}M^{-1}}}
\end{equation*}
and the right hand side converges to $M/4$ when $\alpha \to 2$, hence in the limit we get $c \geq M/4$, which is exactly the same bound that can be obtained in the diffusion case (see \cite{eberle} once again).
\end{example}

\section*{Acknowledgement}
I would like to thank my PhD advisor, Andreas Eberle, for numerous discussions and valuable comments regarding the contents of this paper. I am also grateful to Zdzis\l aw Brze\'{z}niak for some useful remarks and to the anonymous referee for constructive suggestions that helped me to improve the paper. This work was financially supported by the Bonn International Graduate School of Mathematics.

\end{document}